\documentclass[12pt]{amsart}
\usepackage{amssymb,verbatim}

\newtheorem{thm}{Theorem}[section]

\newtheorem{fact}[thm]{Fact}

\newtheorem{cor}[thm]{Corollary}

\newtheorem{claim}{Claim}[thm]

\theoremstyle{remark}

\newtheorem*{remark}{Remark}

\theoremstyle{definition}
\newtheorem{defn}[thm]{Definition}

\newtheorem{q}{Question}
\newtheorem*{qq}{Question}

\def\s{\subseteq}
\def\bks{\setminus}

\def\gorer{\Rightarrow}

\def\goreri{\Leftrightarrow}
\def\notgorer{\not\Rightarrow}
\fboxsep0.1mm\def\sd{\framebox[3.4mm][l]{$\diamondsuit$}\hspace{0.5mm}{}}
\def\ssd{\framebox[2.4mm][l]{$\diamondsuit$}\hspace{0.4mm}{}}

\DeclareMathOperator{\add}{Add}
\DeclareMathOperator{\nacc}{nacc}
\DeclareMathOperator{\tr}{Tr}
\DeclareMathOperator{\cf}{cf}
\DeclareMathOperator{\ma}{MA}
\DeclareMathOperator{\dom}{dom}
\DeclareMathOperator{\otp}{otp} \DeclareMathOperator{\acc}{acc}
\DeclareMathOperator{\sap}{SAP}
\DeclareMathOperator{\zfc}{ZFC}
\DeclareMathOperator{\zf}{ZF}
\DeclareMathOperator{\ad}{AD}
\DeclareMathOperator{\ap}{AP}

\DeclareMathOperator{\ns}{NS} 

\DeclareMathOperator{\ch}{CH}

\DeclareMathOperator{\gch}{GCH}
\DeclareMathOperator{\refl}{Refl}
\DeclareMathOperator{\mm}{\text{\sf MM}}
\DeclareMathOperator{\pfa}{\text{\sf PFA}}

\makeindex

\begin{document}

\title{Jensen's diamond principle and its relatives}

\author{Assaf Rinot}
\address{School of Mathematical Sciences\\Tel Aviv University\\Tel Aviv 69978, Israel}
\email{survey01@rinot.com} \urladdr{http://www.assafrinot.com}


\begin{abstract} We survey some recent results on the validity of Jensen's diamond principle at successor cardinals.
We also discuss weakening of this principle such as club guessing, and anti-diamond principles such as uniformization.

A collection of open problems is included.
\end{abstract}

\maketitle
\tableofcontents
\subsection*{Introduction}
 Cantor's continuum hypothesis has many equivalent formulations in the context of $\zfc$.
One of the standard formulations asserts the existence of an enumeration $\{ A_\alpha\mid \alpha<\omega_1\}$
of the set $\mathcal P(\omega)$. A non-standard, twisted, formulation of $\ch$ is as follows:
\begin{itemize}
\item[$(\exists)_2$] there exists a sequence, $\langle A_\alpha\mid \alpha<\omega_1\rangle$, such that for every subset $Z\s\omega_1$,
there exist two infinite ordinals $\alpha,\beta<\omega_1$ such that $Z\cap\beta=A_\alpha$.
\end{itemize}

By omitting one of the two closing quantifiers in the above statement, we arrive to the following enumeration principle:
\begin{itemize}
\item[$(\exists)_1$] there exists a sequence, $\langle A_\alpha\mid \alpha<\omega_1\rangle$, such that for every subset $Z\s\omega_1$,
there exists an infinite ordinal $\alpha<\omega_1$ such that $Z\cap\alpha=A_\alpha$.
\end{itemize}

Jensen discovered this last principle during his analysis of G\"odel's constructible universe,
and  gave it the name of \emph{diamond}, $\diamondsuit$.
In \cite{jensen}, Jensen proved that $\diamondsuit$ holds in the constructible universe,
and introduced the very first $\diamondsuit$-based construction of a complicated combinatorial object --- a Souslin tree.
Since then, this principle and generalizations of it became very popular among set theorists
who utilized it to settle open problems in fields including topology, measure theory and group theory.

In this paper, we shall be discussing a variety of diamond-like
principles for successor cardinals, including \emph{weak diamond},  \emph{middle diamond},
\emph{club guessing}, \emph{stationary hitting}, and \emph{$\lambda^+$-guessing},
as well as, anti-diamond principles 
including the \emph{uniformization property} and the \emph{saturation of the nonstationary ideal}.

An effort has been put toward including a lot of material, while maintaining an healthy reading flow.
In particular,
this survey cannot cover all known results on this topic.
Let us now briefly describe the content of this survey's sections, and comment on the chosen focus of each section.

\subsection*{Organization of this paper}

In Section 1, Jensen's diamond principles, $\diamondsuit_S$, $\diamondsuit^*_S$, $\diamondsuit^+_S$, are discussed.
We address the question to which stationary sets $S\s\lambda^+$, does $2^\lambda=\lambda^+$ imply $\diamondsuit_S$ and $\diamondsuit^*_S$,
and describe the effect of square principles and reflection principles on diamond.
We discuss a $\gch$-free version of diamond, which is called \emph{stationary hitting},
and a reflection-free version of $\diamondsuit^*_S$, denoted by $\diamondsuit^{\lambda^+}_S$.
In this section, we only deal with the most fundamental variations of diamond, and hence we can outline the whole history.

Section 2 is dedicated to describing part of the set theory generated by Whitehead problem.
We deal with the weak diamond, $\Phi_S$, and the uniformization property.
Here, rather than including all known results in this direction, we decided to focus on presenting the illuminating proofs
of the characterization of weak diamond in cardinal-arithmetic terms, and the failure of instances of  the uniformization property
at successor of singular cardinals.

In Section 3, we go back to the driving force to the study of diamond --- the Souslin hypothesis.
Here, we  focus on aggregating old, as well as, new open problems around the existence of higher souslin trees,
and the existence of particular club guessing sequences.

Section 4 deals with non-saturation of particular ideals --- ideals of the form $\ns_{\lambda^+}\restriction S$.
Here, we describe the interplay between non-saturation, diamond and weak-diamond,
and we focus on presenting the recent results in this line of research.

\subsection*{Notation and conventions}
For ordinals $\alpha<\beta$, we denote by $(\alpha,\beta):=\{\gamma\mid \alpha<\gamma<\beta\}$,
the open interval induced by $\alpha$ and $\beta$.
For a set of ordinals $C$,
we denote by $\acc(C):=\{ \alpha<\sup(C)\mid \sup(C\cap\alpha)=\alpha\}$,
 the set of all accumulation points of $C$.
For a regular uncountable cardinal, $\kappa$, and a subset $S\s\kappa$, let
$$\tr(S):=\min\{\gamma<\kappa\mid \cf(\gamma)>\omega, S\cap\gamma\text{ is stationary in }\gamma\}.$$
We say that $S$ \emph{reflects} iff $\tr(S)\not=\emptyset$, is \emph{non-reflecting} iff $\tr(S)=\emptyset$,
 and \emph{reflects  stationarily often}  iff $\tr(S)$ is stationary.

For cardinals $\kappa<\lambda$, denote $E^{\lambda}_{\kappa}:=\left\{ \alpha<\lambda\mid \cf(\alpha)=\kappa\right\}$,
and $[\lambda]^\kappa:=\{ X\s\lambda\mid |X|=\kappa\}$. $E^{\lambda}_{>\kappa}$ and $[\lambda]^{<\kappa}$ are defined analogously.
Cohen's notion of forcing for adding $\kappa$ many $\lambda$-Cohen sets is denoted by $\add(\lambda,\kappa)$.
To exemplify, the forcing notion for adding a single Cohen real is denoted by $\add(\omega,1)$.

\section{Diamond}
Recall Jensen's notion of diamond in the context of successor cardinals.

\begin{defn}[Jensen, \cite{jensen}]\label{0} For an infinite cardinal $\lambda$ and stationary subset $S\s\lambda^+$:
\begin{itemize}
\item[$\blacktriangleright$]\index{Guessing Principles!$\diamondsuit_S$!Set version} $\diamondsuit_S$ asserts that there exists a sequence $\langle A_\alpha\mid \alpha\in S\rangle$ such that:
\begin{itemize}
\item[$\bullet$] for all $\alpha\in S$, $A_\alpha\s\alpha$;
\item[$\bullet$] if $Z$ is a subset of $\lambda^+$, then the following set is stationary:
$$\{ \alpha\in S\mid Z\cap\alpha=A_\alpha\}.$$
\end{itemize}
\end{itemize}
\end{defn}

Jensen isolated the notion of diamond  from his original construction of an $\aleph_1$-Souslin tree from $V=L$;
in \cite{jensen}, he proved that $\diamondsuit_{\omega_1}$ witnesses the existence of such a tree, and that:

\begin{thm}[Jensen, \cite{jensen}]\label{01} If $V=L$, then $\diamondsuit_S$ holds for every stationary $S\s\lambda^+$
and every infinite cardinal $\lambda$.
\end{thm}

In fact, Jensen established that  $V=L$ entails stronger versions of diamond,
two of which are the following.

\begin{defn}[Jensen, \cite{jensen}]\label{def13} For an infinite cardinal $\lambda$ and stationary subset $S\s\lambda^+$:
\begin{itemize}
\item[$\blacktriangleright$]\index{Guessing Principles!$\diamondsuit^*_S$} $\diamondsuit^*_S$ asserts that there exists a sequence $\langle \mathcal A_\alpha\mid \alpha\in S\rangle$ such that:
\begin{itemize}
\item[$\bullet$] for all $\alpha\in S$, $\mathcal A_\alpha\s\mathcal P(\alpha)$ and $|\mathcal A_\alpha|\le\lambda$;
\item[$\bullet$] if $Z$ is a subset of $\lambda^+$, then the there exists a club $C\s\lambda^+$ such that:
$$C\cap S\s \{ \alpha\in S\mid Z\cap\alpha\in\mathcal A_\alpha\}.$$
\end{itemize}
\item[$\blacktriangleright$]\index{Guessing Principles!$\diamondsuit^+_S$} $\diamondsuit^+_S$ asserts that there exists a sequence $\langle \mathcal A_\alpha\mid \alpha\in S\rangle$ such that:
\begin{itemize}
\item[$\bullet$] for all $\alpha\in S$, $\mathcal A_\alpha\s\mathcal P(\alpha)$ and $|\mathcal A_\alpha|\le\lambda$;
\item[$\bullet$] if $Z$ is a subset of $\lambda^+$, then the there exists a club $C\s\lambda^+$ such that:
$$C\cap S\s \{ \alpha\in S\mid Z\cap\alpha\in\mathcal A_\alpha\ \&\ C\cap\alpha\in\mathcal A_\alpha\}.$$
\end{itemize}
\end{itemize}
\end{defn}

Kunen \cite{kunen} proved that $\diamondsuit_S^*\gorer \diamondsuit_T$ for every stationary $T\s S\s\lambda^+$,
and that  $\diamondsuit_{\lambda^+}$ cannot be introduced by a $\lambda^+$-c.c. notion of forcing.

Since, for a stationary subset $S\s\lambda^+$, $\diamondsuit^+_S\gorer\diamondsuit^*_S\gorer\diamondsuit_S\gorer \diamondsuit_{\lambda^+}\gorer(2^\lambda=\lambda^+)$,
it is natural to study which of these implications may be reversed.

Jensen (see \cite{MR491861}) established the consistency of $\diamondsuit^*_{\omega_1}+\neg\diamondsuit^+_{\omega_1}$,
from the existence of an inaccessible cardinal.
In \cite{rinot09}, it is observed that if $\lambda^{\aleph_0}=\lambda$,
then for every stationary $S\s\lambda^+$, $\diamondsuit^*_S$ is equivalent to $\diamondsuit^+_S$.
Devlin \cite{devlinvariation}, starting with a model of $V\models\gch$, showed that $V^{\add(\lambda^+,\lambda^{++})}\models \neg\diamondsuit^*_{\lambda^+}+\diamondsuit_{\lambda^+}$.\footnote{
For this, he argued that if $G$ is $\add(\lambda^+,1)$-generic over $V$,
then
\begin{enumerate}
\item   $V[G]\models \diamondsuit_S$ for every stationary $S\s\lambda^+$ from $V$, and
\item  every sequence $\langle \mathcal A_\alpha\mid \alpha<\lambda^+\rangle$ that witnesses $\diamondsuit^*_{\lambda^+}$ in $V$,
will cease to witness $\diamondsuit^*_{\lambda^+}$ in $V[G]$.
\end{enumerate}}
Jensen proved that, in general, the implication $\diamondsuit_{\lambda^+}\gorer(2^\lambda=\lambda^+)$, may not be reversed:

\begin{thm}[Jensen, see \cite{devlins}]\label{14} $\ch$ is consistent together with $\neg\diamondsuit_{\omega_1}$.
\end{thm}

On the other hand, Gregory, in a paper that deals with higher Souslin trees, established the following surprising result.

\begin{thm}[Gregory, \cite{gregory}]\label{15} Suppose $\lambda$ is an uncountable cardinal, $2^\lambda=\lambda^+$.

If $\sigma<\lambda$ is an infinite cardinal such that $\lambda^{\sigma}=\lambda$, then $\diamondsuit^*_{E^{\lambda^+}_{\sigma}}$ holds.

In particular, $\gch$ entails $\diamondsuit^*_{E^{\lambda^+}_{<\cf(\lambda)}}$ for any cardinal $\lambda$ of uncountable cofinality.
\end{thm}

Unfortunately, it is impossible to infer $\diamondsuit_{\lambda^+}$ from $\gch$
using Gregory's theorem, in the case that $\lambda>\cf(\lambda)=\omega$.
However, shortly afterwards, this missing case has been settled by Shelah.

\begin{thm}[Shelah, \cite{sh108}]\label{16} Suppose $\lambda$ is a singular cardinal, $2^\lambda=\lambda^+$.

If $\sigma<\lambda$ is an infinite cardinal such that
$\sup\{ \mu^{\sigma}\mid \mu<\lambda\}=\lambda$,
and $\sigma\not=\cf(\lambda)$, then $\diamondsuit^*_{E^{\lambda^+}_\sigma}$ holds.

In particular, $\gch$ entails $\diamondsuit^*_{E^{\lambda^+}_{\not=\cf(\lambda)}}$ for every uncountable cardinal, $\lambda$.
\end{thm}

A closer look at the proof of Theorems \ref{15}, \ref{16} reveals that moreover  $\diamondsuit^+_{E^{\lambda^+}_\sigma}$ may be inferred from the same assumptions,
and, more importantly, that the hypothesis involving $\sigma$ may be weakened to: ``$\sup\{ \cf([\mu]^{\sigma},\supseteq) \mid \mu<\lambda\}=\lambda$''.
However, it was not clear to what extent this weakening indeed witnesses more instances of diamonds.

Then, twenty years after proving Theorem \ref{16}, Shelah established that the above weakening is quite prevalent.
In \cite{sh460}, he proved that the following consequence of $\gch$ follows outright from $\zfc$.

\begin{thm}[Shelah, \cite{sh460}]\label{17} If $\theta$ is an uncountable strong limit cardinal,
then for every cardinal $\lambda\ge\theta$, the set $\left\{ \sigma<\theta\mid  \cf([\lambda]^{\sigma},\supseteq)>\lambda\right\}$
is bounded below $\theta$.

In particular, for every cardinal $\lambda\ge\beth_\omega$, the following are equivalent:
\begin{enumerate}
\item $2^\lambda=\lambda^+$;
\item $\diamondsuit_{\lambda^+}$;
\item $\diamondsuit^*_{E^{\lambda^+}_\sigma}$ for co-boundedly many $\sigma<\beth_\omega$.
\end{enumerate}
\end{thm}

Let $\ch_\lambda$\index{____Ladd@$\ch_\lambda$} denote the assertion that $2^\lambda=\lambda^+$.
By Theorems \ref{14} and \ref{17},  $\ch_\lambda$ does not imply $\diamondsuit_{\lambda^+}$ for $\lambda=\omega$,
but does imply $\diamondsuit_{\lambda^+}$ for every cardinal $\lambda\ge\beth_\omega$.
This left a mysterious gap between $\omega$ and $\beth_\omega$,
which was only known to be closed in the presence of the stronger cardinal arithmetic hypotheses, as in Theorem \ref{15}.

It then took ten additional years until this mysterious gap has been completely closed,
where recently Shelah proved the following striking theorem.

\begin{thm}[Shelah, \cite{sh922}]\label{13} For an uncountable cardinal $\lambda$, and a stationary subset $S\s E^{\lambda^+}_{\not=\cf(\lambda)}$,
the following are equivalent:
\begin{enumerate}
\item $\ch_\lambda$;
\item $\diamondsuit_S$.
\end{enumerate}
\end{thm}

\begin{remark}
In \cite{komjath}, Komj\'ath provides a simplified presentation of Shelah's proof.
Also, in \cite{rinot07} the author presents a considerably shorter proof.\footnote{See the discussion after Theorem \ref{12} below.}
\end{remark}
Having Theorem \ref{13} in hand, we now turn to studying the validity of $\diamondsuit_S$
for sets of the form $S\s E^{\lambda^+}_{\cf(\lambda)}$.
Relativizing Theorem \ref{15} to the
first interesting case, the case $\lambda=\aleph_1$,
we infer that $\gch$ entails $\diamondsuit^*_{E^{\omega_2}_\omega}$.
By Devlin's theorem \cite{devlinvariation}, $\gch\notgorer \diamondsuit^*_{\aleph_2}$,
and consequently, $\gch$ does not imply $\diamondsuit^*_{E^{\aleph_2}_{\omega_1}}$.
Now, what about the unstarred version of diamond?
It turns out that the behavior here is analogous to the one of Theorem \ref{14}.

\begin{thm}[Shelah, see \cite{{kingsteinhorn}}]\label{109} $\gch$ is consistent with $\neg\diamondsuit_{S}$, for $S=E^{\omega_2}_{\omega_1}$.
\end{thm}

The proof of Theorem \ref{109} generalizes to successor of higher regular cardinals,
suggesting that we should focus our attention on successors of singulars. And indeed,
a longstanding, still open, problem is the following question.

\begin{q}[Shelah]\label{q1}\index{ZZZ@Open Problems!Question 01} Is it consistent that
for some \emph{singular} cardinal $\lambda$, $\ch_\lambda$ holds,
while $\diamondsuit_{E^{\lambda^+}_{\cf(\lambda)}}$ fails?
\end{q}

In \cite[$\S3$]{sh186}, Shelah established that a positive answer to the above question
--- in the case that $\lambda$ is a strong limit ---
would entail the failure of weak square,\footnote{\index{Square Principles!$\square^*_\lambda$}The weak square property at $\lambda$, denoted $\square^*_\lambda$,
is the principle $\square_{\lambda,\lambda}$ as in Definition \ref{45}.}
and hence requires large cardinals. More specifically:
\begin{thm}[Shelah, \cite{sh186}]\label{110} Suppose $\lambda$ is a strong limit singular cardinal,
and $\square^*_\lambda$ holds.
If $S\s E^{\lambda^+}_{\cf(\lambda)}$ reflects stationarily often, then $\ch_\lambda\gorer\diamondsuit_S$.
\end{thm}

Applying ideas of the proof of Theorem \ref{13} to the proof Theorem \ref{110},
Zeman established a ``strong limit''-free version of the preceding.
\begin{thm}[Zeman, \cite{zeman}]\label{111} Suppose $\lambda$ is a singular cardinal,
and $\square^*_\lambda$ holds.

If $S\s E^{\lambda^+}_{\cf(\lambda)}$ reflects stationarily often, then $\ch_\lambda\gorer\diamondsuit_S$.
\end{thm}

The curious reader may wonder on the role of the reflection hypothesis in the preceding two theorems;
in  \cite[$\S2$]{sh186}, Shelah established the following counterpart:
\begin{thm}[Shelah, \cite{sh186}]\label{112}  Suppose $\ch_\lambda$ holds for a strong limit singular cardinal, $\lambda$.
If $S\s E^{\lambda^+}_{\cf(\lambda)}$ is a non-reflecting stationary set,
then there exists a notion of forcing $\mathbb P_S$ such that:
\begin{enumerate}
\item $\mathbb P_S$ is $\lambda$-distributive;
\item $\mathbb P_S$ satisfies the $\lambda^{++}$-c.c.;
\item $S$ remains stationary in $V^{\mathbb P_S}$;
\item $V^{\mathbb P_S}\models \neg\diamondsuit_S$.
\end{enumerate}

In particular, it is consistent that $\gch+\square^*_\lambda$ holds, while $\diamondsuit_S$
fails for some non-reflecting stationary set $S\s E^{\lambda^+}_{\cf(\lambda)}$.
\end{thm}

The next definition suggests a way of filtering out the behavior of diamond on non-reflecting sets.

\begin{defn}[\cite{rinot07}] For an infinite cardinal $\lambda$ and stationary subsets $T,S\s\lambda^+$:
\begin{itemize}
\item[$\blacktriangleright$]\index{Guessing Principles!$\diamondsuit^T_S$} $\diamondsuit^T_S$ asserts that there exists a sequence $\langle \mathcal A_\alpha\mid \alpha\in S\rangle$ such that:
\begin{itemize}
\item[$\bullet$] for all $\alpha\in S$, $\mathcal A_\alpha\s\mathcal P(\alpha)$ and $|\mathcal A_\alpha|\le\lambda$;
\item[$\bullet$] if $Z$ is a subset of $\lambda^+$, then the following set is non-stationary:
$$T\cap\tr\{ \alpha\in S\mid Z\cap\alpha\not\in\mathcal A_\alpha\}.$$
\end{itemize}
\end{itemize}
\end{defn}

Notice that by Theorem \ref{16}, $\gch$ entails $\diamondsuit^{\lambda^+}_{\lambda^+}$ for every regular cardinal $\lambda$.
Now, if $\lambda$ is singular, then $\gch$ does not necessarily imply $\diamondsuit^{\lambda^+}_{\lambda^+}$,%
\footnote{Start with a model of $\gch$ and a supercompact cardinal $\kappa$.
Use backward Easton support iteration of length $\kappa+1$, forcing with $\add(\alpha^{+\omega+1},\alpha^{+\omega+2})$
for every inaccessible $\alpha\le\kappa$. Now, work in the extension and let $\lambda:=\kappa^{+\omega}$. Then the $\gch$ holds,
$\kappa$ remains supercompact,
and by Devlin's argument \cite{devlinvariation}, $\diamondsuit^*_{\lambda^+}$ fails. Since $\cf(\lambda)<\kappa<\lambda$,
and $\kappa$ is supercompact, we get that every stationary subset of $E^{\lambda^+}_{\cf(\lambda)}$ reflects, and so it follows from Theorem \ref{113}(2),
that $\diamondsuit^{\lambda^+}_{\lambda^+}$ fails in this model of $\gch$.}
however, if in addition $\square^*_\lambda$ holds,
then $\gch$ does entail $\diamondsuit^{\lambda^+}_{\lambda^+}$,
as the following improvement of theorem \ref{110} shows.
\begin{thm}[\cite{rinot07}]\label{113} For a strong limit singular cardinal, $\lambda$:
\begin{enumerate}
\item if $\square^*_\lambda$ holds, then  $\ch_\lambda\goreri\diamondsuit^{\lambda^+}_{\lambda^+}$;
\item if every stationary subset of $E^{\lambda^+}_{\cf(\lambda)}$ reflects, then $\diamondsuit^{\lambda^+}_{\lambda^+}\goreri\diamondsuit^*_{\lambda^+}$.
\end{enumerate}
\end{thm}

\begin{remark}
An interesting consequence of the preceding theorem is that assuming $\gch$,
for every singular cardinal, $\lambda$,
$\square^*_\lambda$ implies
that in the generic extension by $\add(\lambda^+,1)$, there exists a non-reflecting stationary subset of $\lambda^+$.
This is a reminiscent of the fact that $\square_\lambda$ entails the existence non-reflecting stationary subset of $\lambda^+$.
\end{remark}

Back to Question \ref{q1}, it is natural to study to what extent can the weak square hypothesis in Theorem \ref{111}
be weakened. We now turn to defining the axiom $\sap_\lambda$ and describing its relation to weak square and diamond.

\begin{defn}[\cite{rinot07}]\index{Square Principles!$I[S;\lambda]$} For a singular cardinal $\lambda$ and $S\subset \lambda^+$, consider the ideal $I[S;\lambda]$:
a set $T$ is in $I[S;\lambda]$ iff $T\s\tr(S)$ and there exists
a function $d:[\lambda^+]^2\rightarrow\cf(\lambda)$
such that:
\begin{itemize}
\item $d$ is \emph{subadditive}:  $\alpha<\beta<\gamma<\lambda^+$ implies $d(\alpha,\gamma)\le\max\{ d(\alpha,\beta),d(\beta,\gamma)\}$;
\item $d$ is \emph{normal}: for all $i<\cf(\lambda)$ and  $\beta<\lambda^+$, $|\{ \alpha<\beta\mid d(\alpha,\beta)\le i\}|<\lambda$;
\item \emph{key property}: for some club $C\s\lambda^+$, for every $\gamma\in T\cap C\cap E^{\lambda^+}_{>\cf(\lambda)}$,
there exists a stationary $S_\gamma\s S\cap\gamma$ with $\sup(d``[S_\gamma]^2)<\cf(\lambda)$.
\end{itemize}
\end{defn}

Evidently, if $I[S;\lambda]$ contains a stationary set, then $S$ reflects stationarily often.
The purpose of the next definition is to impose the converse implication.
\begin{defn}[\cite{rinot07}]\index{Square Principles!$\sap_\lambda$} For a singular cardinal $\lambda$, the \emph{stationary approachability property} at $\lambda$,
abbreviated $\sap_\lambda$, asserts that
$I[S;\lambda]$ contains a stationary set for every stationary $S\s E^{\lambda^+}_{\cf(\lambda)}$
that reflects stationarily often.
\end{defn}

Our ideal $I[S;\lambda]$ is a variation of Shelah's approachability ideal $I[\lambda^+]$,
and the axiom $\sap_\lambda$ is a variation of the \emph{approachability property}, $\ap_\lambda$.\footnote{%
For instance, if $\lambda>\cf(\lambda)>\omega$ is a strong limit,
then $I[\lambda^+]=\mathcal P(E^{\lambda^+}_\omega)\cup I[E^{\lambda^+}_{\omega};\lambda]$.
For the definition of $I[\lambda^+]$ and $\ap_\lambda$, see \cite{eisworth}.}
We shall be comparing these two principles later, but let us first compare $\sap_\lambda$ with $\square^*_\lambda$.

In \cite{rinot07}, it is observed that for every singular cardinal $\lambda$, $\square^*_\lambda\gorer\sap_\lambda$,
and moreover, $\square^*_\lambda$ entails the existence of a function, $d:[\lambda^+]^2\rightarrow\cf(\lambda)$,
that serves as a \emph{unified} witness to the fact \emph{for all} $S\s \lambda^+$, $\tr(S)\in I[S;\lambda]$.
Then, starting with a supercompact cardinal, a model is constructed in which (1) $\gch+\sap_{\aleph_\omega}$ holds,
(2) every stationary subset of $E^{\aleph_{\omega+1}}_\omega$ reflects stationarily often,
and (3) for every stationary $S\s E^{\aleph_{\omega+1}}_{\omega}$ and any function $d$
witnessing that $I[S;\aleph_\omega]$ contains a stationary set,
there exists another stationary $S'\s E^{\aleph_{\omega+1}}_{\omega}$
such that this particular $d$ does not witness the fact that $I[S';\aleph_\omega]$ contains a stationary set.
Thus, establishing:
\begin{thm}[\cite{rinot07}]\label{117}
It is relatively consistent with the existence of a supercompact cardinal,
that $\sap_{\aleph_\omega}$ holds, while $\square^*_{\aleph_\omega}$ fails.
\end{thm}

Once it is established that $\sap_\lambda$ is strictly weaker than $\square^*_\lambda$,
the next task would be proving that it is possible to replace $\square^*_\lambda$ in Theorem \ref{111}
with $\sap_\lambda$, while obtaining the same conclusion.
The proof of this fact goes through a certain cardinal-arithmetic-free version of diamond, which we now turn to define.

\begin{defn}[\cite{rinot07}]  For an infinite cardinal $\lambda$ and stationary subsets $T,S\s\lambda^+$,
consider the following two principles:
\begin{itemize}
\item[$\blacktriangleright$]\index{Guessing Principles!Z1@$\clubsuit^-_S$!Set version} $\clubsuit^-_S$ asserts that there exists a sequence $\langle \mathcal A_\alpha\mid \alpha\in S\rangle$ such that:
\begin{itemize}
\item[$\bullet$] for all $\alpha\in S$, $\mathcal A_\alpha\s [\alpha]^{<\lambda}$ and $|\mathcal A_\alpha|\le\lambda$;
\item[$\bullet$] if $Z$ is a \emph{cofinal} subset of $\lambda^+$, then the following set is stationary:
$$\left\{\alpha\in S\mid \exists A\in\mathcal A_\alpha( \sup(Z\cap A)=\alpha)\right\}.$$
\end{itemize}

\item[$\blacktriangleright$]\index{Guessing Principles!ZS@Stationary hitting} $\clubsuit^-_S\restriction T$ asserts that there exists a sequence $\langle \mathcal A_\alpha\mid \alpha\in S\rangle$ such that:
\begin{itemize}
\item[$\bullet$] for all $\alpha\in S$, $\mathcal A_\alpha\s [\alpha]^{<\lambda}$ and $|\mathcal A_\alpha|\le\lambda$;
\item[$\bullet$] if $Z$ is a \emph{stationary} subset of $T$, then the following set is non-empty:
$$\left\{\alpha\in S\mid \exists A\in\mathcal A_\alpha( \sup(Z\cap A)=\alpha)\right\}.$$
\end{itemize}
\end{itemize}
\end{defn}

Notice that $\clubsuit^-_S$ makes sense only in the case that $S\s E^{\lambda^+}_{<\lambda}$.
In \cite{rinot07}, it is established that the \emph{stationary hitting} principle, $\clubsuit^-_S\restriction\lambda^+$,
is equivalent to $\clubsuit^-_S$, and that these equivalent principles are the cardinal-arithmetic-free version of diamond:

\begin{thm}[\cite{rinot07}]\label{12} For an uncountable cardinal $\lambda$, and a stationary subset $S\s E^{\lambda^+}_{<\lambda}$,
the following are equivalent:
\begin{enumerate}
\item $\clubsuit^-_S+\ch_\lambda$;
\item $\diamondsuit_S$.
\end{enumerate}
\end{thm}

It is worth mentioning that the proof of Theorem \ref{12} is surprisingly short,
and when combined with the easy argument that
$\zfc\vdash \clubsuit^-_S$ for every stationary subset $S\s E^{\lambda^+}_{\not=\cf(\lambda)}$,
one obtains a single-page proof of Theorem \ref{13}.

It is also worth mentioning the functional versions of these principles.

\begin{fact}\label{120} Let $\lambda$ denote an infinite cardinal, and $S$ denote a stationary subset of $\lambda^+$; then:
\begin{enumerate}
\item[$\blacktriangleright$]\index{Guessing Principles!$\diamondsuit_S$!Functional version}  $\diamondsuit_S$ is equivalent to the existence of a sequence $\langle g_\alpha\mid \alpha\in S\rangle$
such that:
\begin{itemize}
\item[$\bullet$] for all $\alpha\in S$, $g_\alpha:\alpha\rightarrow\alpha$ is some function;
\item[$\bullet$]  for every function $f:\lambda^+\rightarrow\lambda^+$, the following set is stationary:
$$\{ \alpha\in S\mid f\restriction\alpha=g_\alpha\}.$$
\end{itemize}
\item[$\blacktriangleright$]\index{Guessing Principles!Z1@$\clubsuit^-_S$!Functional version} $\clubsuit^-_S$ is equivalent to the existence of a sequence
$\langle \mathcal G_\alpha\mid \alpha\in S\rangle$ such that:
\begin{itemize}
\item[$\bullet$] for all $\alpha\in S$, $\mathcal G_\alpha\s[\alpha\times\alpha]^{<\lambda}$ and $|\mathcal G_\alpha|\le\lambda$;
\item[$\bullet$]  for every function $f:\lambda^+\rightarrow\lambda^+$,
the following set is stationary:
$$\left\{\alpha\in S\mid \exists G\in\mathcal G_\alpha\ \sup\{\beta<\alpha\mid (\beta,f(\beta))\in G\}=\alpha\right\}.$$
\end{itemize}
\end{enumerate}
\end{fact}

Finally, we are now in a position to formulate a theorem of local nature, from which we derive a global corollary.

\begin{thm}[\cite{rinot07}]\label{18} Suppose $\lambda$ is a singular cardinal,
and  $S\s \lambda^+$ is a stationary set.
If $I[S;\lambda]$ contains a stationary set, then $\clubsuit^-_S$ holds.
\end{thm}
\begin{cor}[\cite{rinot07}]\label{121}Suppose $\sap_\lambda$ holds, for a given singular cardinal, $\lambda$.
Then the following are equivalent:
\begin{enumerate}
\item $\ch_\lambda$;
\item $\diamondsuit_S$ holds for every $S\s \lambda^+$ that reflects stationarily often.
\end{enumerate}
\end{cor}

Thus, the hypothesis $\square^*_\lambda$ from Theorem \ref{111} may indeed be weakened to $\sap_\lambda$.
Having this positive result in mind, one may hope to improve the preceding, proving that $\ch_\lambda\gorer\diamondsuit_S$
for every $S\s\lambda^+$ that reflects stationarily often, without any additional assumptions.
Clearly, this would have settle Question 1 (in the negative!).
However, a recent result by  Gitik and the author shows that diamond \emph{may fail} on a set that reflects
stationarily often, and even on an $(\omega_1+1)$-fat subset of $\aleph_{\omega+1}$:
\begin{thm}[Gitik-Rinot, \cite{gitik-rinot}]\label{124} It is relatively consistent with the existence of a supercompact cardinal that
the $\gch$ holds above $\omega$, while
$\diamondsuit_S$ fails for a stationary set $S\s E^{\aleph_{\omega+1}}_\omega$ such that:
$$\{\gamma<\aleph_{\omega+1}\mid \cf(\gamma)=\omega_1, S\cap\gamma\text{ contains a club}\}\text{ is stationary}.$$
\end{thm}

In fact, the above theorem is just one application of the following general, $\zfc$ result.

\begin{thm}[Gitik-Rinot, \cite{gitik-rinot}]
\label{19} Suppose $\ch_\lambda$ holds for a strong limit singular cardinal, $\lambda$.
Then there exists a notion of forcing $\mathbb P$, satisfying:
\begin{enumerate}
\item $\mathbb P$ is $\lambda^+$-directed closed;
\item $\mathbb P$ has the $\lambda^{++}$-c.c.;
\item $|\mathbb P|=\lambda^{++}$;
\item in $V^{\mathbb P}$, $\diamondsuit_S$ fails for some stationary $S\s E^{\lambda^+}_{\cf(\lambda)}$.
\end{enumerate}
\end{thm}

Note that unlike Theorem \ref{113}, here the stationary set on which diamond fails, is a generic one.

Utilizing the forcing notion from Theorem \ref{19}, Gitik and the author were able to show that Corollary \ref{121}
is optimal:  in \cite{gitik-rinot}, it is proved that replacing the $\sap_\lambda$ hypothesis in Corollary \ref{121} with $\ap_\lambda$,
or with the existence of a \emph{better scale} for $\lambda$, or even with the
existence of a \emph{very good scale} for $\lambda$, is impossible,
in the sense that these alternative hypotheses do not entail diamond on all reflecting stationary sets.\footnote{The existence of a better scale at $\lambda$,
as well as the approachability property at $\lambda$, are well-known consequences of $\square^*_\lambda$. For definitions and proofs, see \cite{eisworth}.}
In particular:

\begin{thm}[Gitik-Rinot, \cite{gitik-rinot}]\label{122} It is relatively consistent with the existence of a
supercompact cardinal that $\ap_{\aleph_\omega}$ holds, while $\sap_{\aleph_\omega}$ fails.
\end{thm}

Moreover, in the model from Theorem \ref{122}, every stationary subset of $E^{\aleph_{\omega+1}}_\omega$ reflects.
Recalling that $\ap_{\aleph_\omega}$ holds whenever every stationary subset of $\aleph_{\omega+1}$ reflects,
we now arrive to the following nice question.

\begin{q}\label{q2}\index{ZZZ@Open Problems!Question 02} Is it consistent that every stationary subset of $\aleph_{\omega+1}$ reflects, while $\sap_{\aleph_\omega}$ fails to hold?
\end{q}

To summarize the effect of  square-like principles on diamond,
we now state a corollary.
Let $\refl_\lambda$
denote the assertion that every stationary subset of $E^{\lambda^+}_{\cf(\lambda)}$ reflects stationarily often. Then:
\begin{cor} For a singular cardinal, $\lambda$:
\begin{enumerate}
\item $\gch+\square_\lambda^*\not\gorer\diamondsuit^*_{\lambda^+}$;
\item $\gch+\refl_\lambda+\square_\lambda^*\gorer\diamondsuit^*_{\lambda^+}$;
\item $\gch+\refl_\lambda+\sap_\lambda\not\gorer\diamondsuit^*_{\lambda^+}$;
\item $\gch+\refl_\lambda+\sap_\lambda\gorer\diamondsuit_{S}$ for every stationary $S\s\lambda^+$;
\item $\gch+\refl_\lambda+\ap_\lambda\not\gorer\diamondsuit_{S}$ for every stationary $S\s\lambda^+$.
\end{enumerate}
\end{cor}
\begin{proof} (1) By Theorem \ref{112}.
(2) By Theorem \ref{113}.
(3) By the proof of Theorem \ref{117} in \cite{rinot07}.
(4) By Corollary \ref{122}.
(5) By the proof of Theorem \ref{122} in \cite{gitik-rinot}.
\end{proof}

The combination of Theorems \ref{12} and \ref{18} motivates the study of the ideal $I[S;\lambda]$.
For instance, a positive answer to the next question would supply an answer to Question \ref{q1}.
\begin{q}\label{q4}\index{ZZZ@Open Problems!Question 03} Must $I[E^{\lambda^+}_{\cf(\lambda)};\lambda]$ contain a stationary set for every singular cardinal $\lambda$?
\end{q}

One of the ways of attacking the above question involves the following reflection principles.
\begin{defn}[\cite{rinot07}] Assume $\theta>\kappa$ are regular uncountable cardinals.

\index{Reflection Principles!$R_1(\theta,\kappa)$}$R_1(\theta,\kappa)$ asserts that for every function $f:E^\theta_{<\kappa}\rightarrow\kappa$,
there exists some $j<\kappa$ such that $\{ \delta\in E^\theta_\kappa\mid f^{-1}[j]\cap\delta\text{ is stationary}\}$ is stationary in $\theta$.

\index{Reflection Principles!$R_2(\theta,\kappa)$}$R_2(\theta,\kappa)$ asserts that for every function $f:E^\theta_{<\kappa}\rightarrow\kappa$,
there exists some $j<\kappa$ such that $\{ \delta\in E^\theta_\kappa\mid f^{-1}[j]\cap\delta\text{ is non-stationary}\}$ is non-stationary.
\end{defn}

It is not hard to see that $R_2(\theta,\kappa)\Rightarrow R_1(\theta,\kappa)$,
and that $\mm$ implies $R_1(\aleph_2,\aleph_1)+\neg R_2(\aleph_2,\aleph_1)$.
In \cite{rinot07}, a fact from \textit{pcf} theory is utilized to prove:
\begin{thm}[\cite{rinot07}] Suppose $\lambda>\cf(\lambda)=\kappa>\omega$ are given cardinals.

The ideal $I[E^{\lambda^+}_{\cf(\lambda)};\lambda]$ contains a stationary set
whenever the following set is non-empty:
$$\{ \theta<\lambda\mid R_1(\theta,\kappa)\text{ holds}\}.$$
\end{thm}

As a corollary, one gets a surprising result stating that a local instance of reflection
affects the validity of diamond on a proper class of cardinals.
\begin{cor}[implicit in \cite{sh922}]
Suppose $\kappa$ is the successor of a cardinal $\kappa^-$,
and that every stationary subset of $E^{\kappa^+}_{\kappa^-}$ reflects.

Then, $\ch_\lambda\goreri\diamondsuit_{E^{\lambda^+}_{\cf(\lambda)}}$ for every singular cardinal $\lambda$ of cofinality $\kappa$.
\end{cor}

As the reader may expect, the principle $R_2$ yields a stronger consequence.

\begin{thm}[\cite{rinot07}]\label{130} Suppose $\theta>\kappa$ are cardinals such that $R_2(\theta,\kappa)$ holds. Then:
\begin{enumerate}
\item For every singular cardinal $\lambda$ of cofinality $\kappa$, and every $S\s\lambda^+$,
we have $$\tr(S)\cap E^{\lambda^+}_\theta\in I[S;\lambda].$$
\item if $\lambda$ is a strong limit singular cardinal of cofinality $\kappa$, then $\ch_\lambda\goreri\diamondsuit_{\lambda^+}^{E^{\lambda^+}_\theta}$.
\end{enumerate}
\end{thm}

Unfortunately, there is no hope to settle Question \ref{q4} using these reflection principles, as they are independent
of \textsf{ZFC}:
by a theorem of Harrington and Shelah \cite{sh99}, $R_1(\aleph_2,\aleph_1)$ is equiconsistent with the existence of a Mahlo cardinal,
whereas, by a theorem of Magidor \cite{magidor}, $R_2(\aleph_2,\aleph_1)$ is consistent modulo the existence of a weakly-compact cardinal.
An alternative sufficient condition for $I[S;\lambda]$ to contain a stationary set will be described in  Section \ref{sectionaturation} (See Fact \ref{131} below).

\newpage
\section{Weak Diamond and the Uniformization Property}\label{sectionwd}

Suppose that $G$ and $H$ are abelian groups and $\pi:H\rightarrow G$ is a given epimorphism.
We say that $\pi$ \emph{splits} iff there exists an homomorphism $\phi:G\rightarrow H$ such that $\pi\circ \phi$ is the identity function on $G$.
An abelian group $G$ is \emph{free} iff every epimorphism onto $G$, splits.

\emph{Whitehead problem} reads as follows.
\begin{qq} Suppose that $G$ is an abelian group such that
every epimorphism $\pi$ onto $G$ with the property that $\ker(\pi)\simeq\mathbb Z$ --- splits;\footnote{Here, $\mathbb Z$ stands for the usual additive group of integers.}

Must $G$ be a free abelian group?
\end{qq}

Thus, the question is whether to decide the freeness of an abelian group,
it suffices to verify that only a particular, narrow, class of epimorphism splits.
Stein \cite{stein} solved Whitehead problem in the affirmative in the case that $G$
is a countable abelian group. Then, in a result that was completely unexpected,
Shelah \cite{sh44} proved that Whitehead problem, restricted to groups of size $\omega_1$, is independent of $\zfc$.
Roughly speaking, by generalizing Stein's proof,
substituting a counting-based diagonalization argument with a guessing-based diagonalization  argument,
Shelah proved that if $\diamondsuit_S$ holds for every stationary $S\s\omega_1$,
then every abelian group of size $\omega_1$ with the above property is indeed free.
On the other hand, he proved that if $\ma_{\omega_1}$ holds, then there exists a counterexample of size $\omega_1$.

Since $\ch$ holds in the first model, and fails in the other, it was natural to ask whether
the existence of a counterexample to Whitehead problem is consistent together with $\ch$.
This led Shelah to introducing the \emph{uniformization property}.

\begin{defn}[Shelah, \cite{sh64}]\label{30001} Suppose that $S$ is a stationary subset of a successor cardinal, $\lambda^+$.
\begin{itemize}
\item[$\bullet$]\index{____Ladder system@Ladder system} A \emph{ladder system} on $S$ is a sequence of sets of ordinals, $\langle L_\alpha\mid \alpha\in S\rangle$,
such that $\sup(L_\alpha)=\alpha$  and $\otp(L_\alpha)=\cf(\alpha)$ for all $\alpha\in S$;

\item[$\bullet$]
\index{Anti-$\diamondsuit$ Principles!The uniformization property}
A ladder system $\langle L_\alpha\mid \alpha\in S\rangle$ is said to have the \emph{uniformization property}
iff whenever $\langle f_\alpha:L_\alpha\rightarrow 2\mid \alpha\in S\rangle$ is a given sequence of local functions,
then there exists a global function $f:\lambda^+\rightarrow 2$ such that $f_\alpha\s ^* f$ for all limit $\alpha\in S$.
That is, $\sup\{\beta\in L_\alpha\mid f_\alpha(\beta)\not=f(\beta)\}<\alpha$ for all limit $\alpha\in S$.
\end{itemize}
\end{defn}

\begin{thm}[Shelah, \cite{sh98}; see also \cite{sh505}]\label{23} The following are equivalent:
\begin{itemize}
\item there exists a counterexample of size $\omega_1$ to Whitehead problem;
\item there exists a stationary $S\s\omega_1$, and a ladder system on $S$
that has the uniformization property.
\end{itemize}
\end{thm}

Devlin and Shelah proved   \cite{sh65} that
if $\ma_{\omega_1}$ holds, then every stationary $S\s\omega_1$
and every ladder system on $S$,
has the uniformization property. On the other hand,
it is not hard to see that if $\diamondsuit_S$ holds, then no ladder system on $S$
has the uniformization property  (See Fact \ref{l25}, below).
Note that altogether, this gives an alternative proof to the independence result from \cite{sh44}.

Recalling that $\neg\diamondsuit_{\omega_1}$ is consistent with $\ch$ (See Theorem \ref{14}),
it seemed reasonable to suspect that $\ch$ is moreover consistent with the existence of a ladder system on $\omega_1$
that has the uniformization property.
Such a model would also show that the existence of a counterexample to Whitehead problem is indeed consistent together with $\ch$,
settling Shelah's question.

However, a surprising theorem of Devlin states that $\ch$ implies that no ladder system on $\omega_1$ has the uniformization property.
Then, in a joint paper with Shelah, the essence of Devlin's proof has been isolated,
and a weakening of diamond which is strong enough to rule out uniformization has been introduced.

\begin{defn}[Devlin-Shelah, \cite{sh65}]  For an infinite cardinal $\lambda$ and a stationary set $S\s\lambda^+$,
consider the principle of \emph{weak diamond}.

\begin{itemize}
\item[$\blacktriangleright$]\index{Guessing Principles!$\Phi_S$}   $\Phi_S$ asserts that
for every function $F:{}^{<\lambda^+}2\rightarrow 2$, there exists a function $g:\lambda^+\rightarrow 2$,
such that for all $f:\lambda^+\rightarrow 2$, the following set is stationary:
$$\{ \alpha\in S\mid F(f\restriction\alpha)=g(\alpha)\}.$$
\end{itemize}
\end{defn}

Note that by Fact \ref{120},  $\diamondsuit_S\gorer\Phi_S$. The difference between these principles is as follows.
In diamond, for each function $f$, we would like to guess $f\restriction\alpha$,
while in weak diamond, we only aim at guessing the value of $F(f\restriction\alpha)$, i.e.,
whether $f\restriction\alpha$ satisfies a certain property --- is it black or white.
A reader who is still dissatisfied with the definition of weak diamond,
may prefer one of its alternative formulations.
\begin{fact}[folklore]
For an infinite cardinal $\lambda$ and a stationary set $S\s\lambda^+$, the following principles are equivalent:
\begin{itemize}
\item[$\blacktriangleright$] $\Phi_S$;
\item[$\blacktriangleright$] for every function $F:{}^{<\lambda^+}\lambda^+\rightarrow 2$, there exists a function $g:S\rightarrow 2$,
such that for all $f:\lambda^+\rightarrow \lambda^+$, the following set is stationary:
$$\{ \alpha\in S\mid F(f\restriction\alpha)=g(\alpha)\}.$$

\item[$\blacktriangleright$] for every sequence of functions $\langle F_\alpha:\mathcal P(\alpha)\rightarrow 2\mid \alpha\in S\rangle$,
there exists a function $g:S\rightarrow 2$, such that for every subset $X\s\lambda^+$, the following set is stationary:
$$\{ \alpha\in S\mid F_\alpha(X\cap\alpha)=g(\alpha)\}.$$

\end{itemize}
\end{fact}

Back to uniformization, we have:

\begin{fact}[Devlin-Shelah, \cite{sh65}]\label{l25} For every stationary set $S$, $\Phi_S$ (and hence $\diamondsuit_S$)
entails that no ladder system $\langle L_\alpha\mid \alpha\in S\rangle$ has the uniformization property.
\end{fact}
\begin{proof}[Proof (sketch)]
For all $\alpha\in S$ and $i<2$, let $c^i_\alpha:L_\alpha\rightarrow\{i\}$ denote the constant function.
Pick a function $F:{}^{<\lambda^+}2\rightarrow 2$ such that for all $\alpha\in S$ and $i<2$,
if $f:\alpha\rightarrow 2$ and  $c^i_\alpha\s^* f$, then $F(f)=i$.
Now, let $g:\lambda^+\rightarrow 2$ be given by applying $\Phi_S$ to $F$.
Then, letting $f_\alpha:=c_\alpha^{1-g(\alpha)}$ for all $\alpha\in S$,
the sequence $\langle f_\alpha\mid \alpha\in S\rangle$ cannot be uniformized.
\end{proof}

Before we turn to showing that $\ch\gorer\Phi_{\omega_1}$, let us mention
that since $\Phi_{\lambda^+}$ deals with two-valued functions, its negation is an interesting statement of its own right:

\begin{fact}\label{34} Suppose that $\Phi_{\lambda^+}$ fails for a given infinite cardinal, $\lambda$.

Then there exists a function $F:{}^{<\lambda^+}({}^\lambda2)\rightarrow{}^\lambda2$
such that for every  $g:\lambda^+\rightarrow{}^\lambda2$, there exists a function $f:\lambda^+\rightarrow{}^\lambda2$,
 for which the following set contains a club:
$$\{\alpha<\lambda^+\mid F(f\restriction\alpha)=g(\alpha)\}.$$
\end{fact}

Roughly speaking, the above states that there exists a decipher, $F$, such that for every function $g$,
there exists a function $f$ that $F$-ciphers the value of $g(\alpha)$ as $f\restriction\alpha$.

Since the (easy) proof of the preceding utilizes the fact that weak diamond deals with two-valued functions,
it is worth mentioning that Shelah also studied generalization involving more colors. For instance,
in \cite{sh638}, Shelah gets weak diamond for more colors provided that $\ns_{\omega_1}$ is saturated (and $\Phi_{\omega_1}$ holds).\footnote{For the definition of ``$\ns_{\omega_1}$ is saturated'' see Definition \ref{40} below.}

We now turn to showing that  $\ch\gorer\Phi_{\omega_1}$. In fact,
the next theorem shows that weak diamond is a cardinal arithmetic statement in disguise.
The proof given here is somewhat lengthier than other available proofs,
but, the value of this proof is that its structure  allows the reader to first
neglect the technical details (by skipping the proofs of Claims \ref{362}, \ref{361}),
while still obtaining a good understanding of the key ideas.

\begin{thm}[Devlin-Shelah, \cite{sh65}]\label{311} For every cardinal $\lambda$, $\Phi_{\lambda^+}\Leftrightarrow 2^{\lambda}<2^{\lambda^+}$.
\end{thm}
\begin{proof} $\Rightarrow$ Assume $\Phi_{\lambda^+}$.
Given an arbitrary function $\psi:{}^{\lambda^+}2\rightarrow{}^\lambda2$,
we now define a function $F:{}^{<\lambda^+}2\rightarrow 2$
such that by appealing to $\Phi_{\lambda^+}$ with $F$, we can show that $\psi$ is not injective.

Given $f\in {}^{<\lambda^+}2$, we let $F(f):=0$ iff there exists a function $h\in{}^{\lambda^+}2$
such that  $h(\dom(f))=0$ and $f\s\psi(h)\cup(h\restriction[\lambda,\lambda^+))$.

Let $g:\lambda^+\rightarrow 2$ be the oracle given by $\Phi_{\lambda^+}$ when applied to $F$,
and let $h:\lambda^+\rightarrow 2$ be the function satisfying $h(\alpha)=1-g(\alpha)$ for all $\alpha<\lambda^+$.

Put $f:=\psi(h)\cup (h\restriction[\lambda,\lambda^+))$. Since $f\in{}^{\lambda^+}2$,
let us pick some  $\alpha<\lambda^+$  with $\alpha>\lambda$ such that $F(f\restriction\alpha)=g(\alpha)$.
Since $f\restriction\alpha\s \psi(h)\cup(h\restriction[\lambda,\lambda^+))$,
the definition of $F$ implies that $F(f\restriction\alpha)=0$ whenever $h(\alpha)=0$.
However, $F(f\restriction\alpha)=g(\alpha)\not=h(\alpha)$, and hence $h(\alpha)=1$.
Since,  $F(f\restriction\alpha)=g(\alpha)=0$, let us  pick a function $h'$ such that
$h'(\alpha)=0$ and $f\restriction\alpha\s\psi(h')\cup(h'\restriction[\lambda,\lambda^+))$.
By definition of $f$, we get that $\psi(h)=f\restriction\lambda=\psi(h')$.
By $g(\alpha)=0$, we also know that $h(\alpha)=1\not=h'(\alpha)$, and hence $h\not=h'$, while $\psi(h)=\psi(h')$.

$\Leftarrow$
Given a function $H:{}^{<\lambda^+}({}^\lambda2)\rightarrow{}^{<\lambda^+}({}^\lambda2)$,
let us say that a sequence $\langle (f_n,D_n)\mid n<\omega\rangle$ is an  $H$-\emph{prospective sequence}
iff:
\begin{enumerate}
\item $\{ D_n\mid n<\omega\}$ is a decreasing chain of club subsets of $\lambda^+$;
\item for all $n<\omega$, $f_n$ is a function from $\lambda^+$ to ${}^\lambda2$;
\item for all $n<\omega$ and $\alpha\in D_{n+1}$, the following holds:
$$H(f_{n+1}\restriction\alpha)=f_n\restriction\min(D_n\bks\alpha+1).$$
\end{enumerate}

Note that the intuitive meaning of the third item is that
there exists $\beta>\alpha$ such that the content of $f_n\restriction\beta$
is coded by $f_{n+1}\restriction\alpha$.

\begin{claim}\label{362} Assume $\neg\Phi_{\lambda^+}$.

Then there exists a function $H:{}^{<\lambda^+}({}^\lambda2)\rightarrow{}^{<\lambda^+}({}^\lambda2)$
such that for every function $f:\lambda^+\rightarrow{}^\lambda2$, there exists an $H$-prospective sequence
$\langle (f_n,D_n)\mid n<\omega\rangle$ with $f_0=f$.
\end{claim}
\begin{proof}
 Fix $F$ as in Fact \ref{34}, and fix a bijection $\varphi:{}^\lambda2\rightarrow{}^{<\lambda^+}(^\lambda2)$.
Put $H:=\varphi\circ F$. Now, given $f:\lambda^+\rightarrow{}^\lambda2$,
we define the $H$-prospective sequence by recursion on $n<\omega$. Start with $f_0:=f$ and $D_0:=\lambda^+$.
Suppose $n<\omega$ and $f_n$ and $D_n$ are defined.
Define a function  $g:\lambda^+\rightarrow {}^{\lambda}2$ by letting for all $\alpha<\lambda^+$:
$$g(\alpha):=\varphi^{-1}(f_n\restriction\min(D_n\bks\alpha+1)).$$
By properties of $F$, there exists a function $f_{n+1}$ and a club $D_{n+1}\s D_n$ such
that for all $\alpha\in D_{n+1}$, we have $F(f_{n+1}\restriction\alpha)=g(\alpha)$.
In particular, $$H(f_{n+1}\restriction\alpha)=(\varphi\circ F)(f_{n+1}\restriction\alpha)=(\varphi\circ g)(\alpha)=f_n\restriction\min(D_n\bks\alpha+1).\qedhere$$

\end{proof}

\begin{claim}\label{361}
Given a function $H:{}^{<\lambda^+}({}^\lambda2)\rightarrow{}^{<\lambda^+}({}^\lambda2)$,
there exists a function $H^*:{}^\omega(^{<\lambda^+}(^\lambda2))\rightarrow{}^\omega(^{<\lambda^+}(^\lambda2))$ with the following stepping-up property.

For every $H$-prospective sequence, $\langle (f_n,D_n)\mid n<\omega\rangle$, and every $\alpha\in\bigcap_{n<\omega}D_n$,
there exists some $\alpha^*<\lambda^+$, such that:
\begin{enumerate}
\item $\alpha^*>\alpha$;
\item $\alpha^*\in \bigcap_{n<\omega}D_n$;
\item $H^*(\langle f_n\restriction\alpha\mid n<\omega\rangle)=\langle f_n\restriction\alpha^*\mid n<\omega\rangle$.
\end{enumerate}

\end{claim}
\begin{proof}
Given  $H$, we define functions $H^m:{}^\omega(^{<\lambda^+}(^\lambda2))\rightarrow{}^\omega(^{<\lambda^+}(^\lambda2))$
by recursion on $m<\omega$. For all $\sigma:\omega\rightarrow{}^{<\lambda^+}(^\lambda2)$, let:
$$H^0(\sigma):=\sigma,$$
and whenever $m<\omega$ is such that $H^m$ is defined,  let:
$$H^{m+1}(\sigma):=\langle H(H^m(\sigma)(n+1)) \mid n<\omega\rangle.$$
Finally, define $H^*$ by letting for all $\sigma:\omega\rightarrow{}^{<\lambda^+}(^\lambda2)$:
$$H^*(\sigma):=\langle \bigcup_{m<\omega}H^m(\sigma)(n)\mid n<\omega\rangle.$$

To see that $H^*$ works, fix an $H$-prospective sequence, $\langle (f_n,D_n)\mid n<\omega\rangle$, and some $\alpha\in\bigcap_{n<\omega}D_n$.
Define $\langle \langle \alpha^m_n\mid n<\omega\rangle\mid m<\omega\rangle$ by letting $\alpha_n^0:=\alpha$
for all $n<\omega$. Then, given $m<\omega$, for all $n<\omega$, let:
$$\alpha_n^{m+1}:=\min(D_n\bks\alpha^m_{n+1}+1).$$

(1) Put $\alpha^*:=\sup_{m<\omega}\alpha^m_0$. Then $\alpha^*\ge\alpha^1_0>\alpha^0_1=\alpha$.

(2) If $n<\omega$, then $D_n\supseteq D_{n+1}$, and hence $\alpha^{m+1}_{n+1}\ge \alpha^{m+1}_n>\alpha^m_{n+1}$ for all $m<\omega$.
This shows that $\sup_{m<\omega}\alpha^m_n=\sup_{m<\omega}\alpha^m_{n+1}$ for all $n<\omega$.

For $n<\omega$, since $\langle \alpha^m_n\mid m<\omega\rangle$
is a strictly increasing sequence of ordinals from $D_n$ that converges to $\alpha^*$, we get that $\alpha^*\in D_n$.

(3) Let us prove by induction that for all $m<\omega$:
$$H^m(\langle f_n\restriction\alpha\mid n<\omega\rangle)=\langle f_n\restriction \alpha^m_n\mid n<\omega\rangle.$$

\underline{Induction Base}: Trivial.

\underline{Induction Step:} Suppose $m<\omega$ is such that:
$$(\star)\qquad H^m(\langle f_n\restriction\alpha\mid n<\omega\rangle)=\langle f_n\restriction \alpha^m_n\mid n<\omega\rangle,$$
and let us show that:
$$H^{m+1}(\langle f_n\restriction\alpha\mid n<\omega\rangle)=\langle f_n\restriction \alpha^{m+1}_n\mid n<\omega\rangle.$$
By definition of $H^{m+1}$ and equation ($\star$), this amounts to showing that:
$$\langle H(f_{n+1}\restriction \alpha^m_{n+1})\mid n<\omega\rangle=\langle f_n\restriction \alpha^{m+1}_n\mid n<\omega\rangle.$$
Fix $n<\omega$. Recalling the definition of $\alpha^{m+1}_n$, we see that we need to prove that $H(f_{n+1}\restriction\alpha^m_{n+1})=f_n\restriction \min(D_n\bks\alpha^{m}_{n+1}+1)$.
But this follows immediately from the facts that $\alpha^m_{n+1}\in D_{n+1}$, and that $\langle (f_n,D_n)\mid n<\omega\rangle$ is an  $H$-prospective sequence.

Thus, it has been established that:
$$H^*(\langle f_n\restriction\alpha\mid n<\omega\rangle)=\langle f_n\restriction\bigcup_{m<\omega}\alpha^m_n\mid n<\omega\rangle=\langle f_n\restriction\alpha^*\mid n<\omega\rangle.\qedhere$$
\end{proof}

Now, assume $\neg\Phi_{\lambda^+}$, and let us prove that $2^{\lambda^+}=2^{\lambda}$
by introducing an injection of the form $\psi:{}^{\lambda^+}(^\lambda2)\rightarrow{}^\omega(^{<\lambda^+}2)$.
Fix $H$ as in Claim \ref{362}, and let $H^*$ be given by Claim \ref{361} when applied to this fixed function, $H$.

$\blacktriangleright$ Given a function $f:\lambda^+\rightarrow{}^\lambda2$, we pick an $H$-prospective sequence $\langle (f_n,D_n)\mid n<\omega\rangle$
with $f_0=f$
and let $\psi(f):=\langle f_n\restriction\alpha\mid n<\omega\rangle$ for $\alpha:=\min(\bigcap_{n<\omega}D_n)$.

To see that $\psi$ is injective, we now define a function $\varphi:{}^\omega(^{<\lambda^+}2)\rightarrow{}^{\le\lambda^+}(^\lambda2)$
such that $\varphi\circ\psi$ is the identity function.

$\blacktriangleright$ Given a sequence $\sigma:\omega\rightarrow{}^{<\lambda^+}2$,
we first define an auxiliary sequence $\langle \sigma_\tau\mid \tau\le\lambda^+\rangle$
by recursion on $\tau$. Let $\sigma_0:=\sigma$, $\sigma_{\tau+1}:=H^*(\sigma_\tau)$,
and $\sigma_{\tau}(n):=\bigcup_{\eta<\tau}\sigma_{\eta}(n)$ for  limit $\tau\le\lambda^+$ and $n<\omega$.
Finally, let $\varphi(\sigma):=\sigma_{\lambda^+}(0)$.

\begin{claim} $\varphi(\psi(f))=f$ for every $f:\lambda^+\rightarrow{}^\lambda2$.
\end{claim}
\begin{proof} Fix $f:\lambda^+\rightarrow{}^\lambda2$ and let $\sigma:=\psi(f)$. By definition of $\psi$,
$\sigma=\langle f_n\restriction\alpha\mid n<\omega\rangle$ for some
 $H$-prospective sequence $\langle (f_n,D_n)\mid n<\omega\rangle$ and $\alpha\in\bigcap_{n<\omega}D_n$.
It then follows from the choice of $H^*$, that there exists a strictly increasing sequence,
 $\langle \alpha_\tau\mid \tau<\lambda^+\rangle$, of ordinals from $\bigcap_{n<\omega}D_n$,
 such that $\sigma_\tau:=\langle f_n\restriction\alpha_\tau\mid n<\omega\rangle$
 for all $\tau<\lambda^+$, and then $\varphi(\psi(f))=\varphi(\sigma)=\sigma_{\lambda^+}(0)=f_0\restriction\lambda^+=f$.
\end{proof}
This completes the proof.
\end{proof}

Evidently, Devlin's pioneering theorem that $\ch$ excludes the existence of
a ladder system on $\omega_1$ with the uniformization property now follows from Fact \ref{l25} and Theorem \ref{311}.
It is interesting to note that if one considers the notion of \emph{weak uniformization},
in which the conclusion of Definition \ref{30001} is weakened from $\sup\{\beta\in L_\alpha\mid f_\alpha(\beta)\not=f(\beta)\}<\alpha$
to $\sup\{\beta\in L_\alpha\mid f_\alpha(\beta)=f(\beta)\}=\alpha$,
then we end up with an example of an anti-$\diamondsuit_S$ principle, which is not an anti-$\Phi_S$ principle:

\begin{thm}[Devlin, see \cite{sh81}]\label{36} It is consistent with $\gch$ (and hence with $\Phi_{\omega_1}$)
that  every ladder system on every stationary subset of $\omega_1$ has the \underline{weak} uniformization property.
\end{thm}

Back to  Whitehead problem, Shelah eventually established the consistency
of $\ch$ together with the existence of a counterexample:

\begin{thm}[Shelah, \cite{sh64}]\label{31} It is consistent with $\gch+\diamondsuit_{\omega_1}$ that there exists a stationary, co-stationary, set $S\s\omega_1$
such that any ladder system on $S$ has the uniformization property.
\end{thm}

It is worth mentioning that Shelah's model was also the first example of a model in which  $\diamondsuit_{\omega_1}$
holds, while for some stationary subset $S\s\omega_1$, $\diamondsuit_S$ fails .

We now turn to dealing with the uniformization property for successor of uncountable cardinals.
By Theorem \ref{13} and Fact \ref{l25}, there is no hope for getting a model of $\gch$ in which a subset of $E^{\lambda^+}_{\not=\cf(\lambda)}$
carries a ladder system that has the uniformization property, so let us focus on sets of the critical cofinality.
The first case that needs to be considered is $E^{\omega_2}_{\omega_1}$,
and the full content of Theorem \ref{109} is now revealed.

\begin{thm}[Shelah, \cite{kingsteinhorn},\cite{sh80}]\label{39} It is consistent with $\gch$ that
there exists a ladder system on $E^{\omega_2}_{\omega_1}$ with the uniformization property.
\end{thm}

Knowing that $2^{\aleph_1}=\aleph_2$ implies $\diamondsuit_{E^{\omega_2}_\omega}$ but not $\diamondsuit_{E^{\omega_2}_{\omega_1}}$,
and that $2^{\aleph_1}<2^{\aleph_2}$ implies $\Phi_{\omega_2}$ but not $\Phi_{E^{\omega_2}_{\omega_1}}$,
one may hope to prove that $2^{\aleph_1}<2^{\aleph_2}$ moreover implies $\Phi_{E^{\omega_2}_{\omega}}$.
However, a consistent counterexample to this conjecture is provided in \cite{sh208}.

Note that Theorem \ref{39} states that there exists a particular ladder system on $E^{\omega_2}_{\omega_1}$
with the uniformization property, rather than stating that all ladder systems on $E^{\omega_2}_{\omega_1}$
have this property.\footnote{Compare with the fact that $\ma_{\omega_1}$ entails that \emph{every} ladder system on $\omega_1$ has the uniformization property.}
To see that Theorem \ref{39} is indeed optimal, consider the following theorem.

\begin{thm}[Shelah, \cite{shpif}]\label{43}
Suppose that $\lambda$ is a regular cardinal of the form $2^\theta$ for some cardinal $\theta$,
and that $\langle L_\alpha\mid \alpha\in E^{\lambda^+}_{\lambda}\rangle$ is a given ladder system.

If, moreover, $L_\alpha$ is a club subset of $\alpha$ for all $\alpha\in E^{\lambda^+}_{\lambda}$,
and $2^{<\lambda}=\lambda$,
then there exists a coloring  $\langle f_\alpha:L_\alpha\rightarrow 2\mid \alpha\in E^{\lambda^+}_{\lambda}\rangle$ such that for every function $f:\lambda^+\rightarrow 2$,
the following set is stationary:
$$\{ \alpha\in E^{\lambda^+}_\lambda\mid \{ \beta\in L_\alpha\mid f_\alpha(\beta)\not=f(\beta)\}\text{ is stationary in }\alpha\}.$$

In particular, $\ch$
entails the existence of a ladder system on $E^{\omega_2}_{\omega_1}$ that does not enjoy the uniformization property.
\end{thm}

The proof of Theorem \ref{39} generalizes to successor of higher regular cardinals,
showing that there may exist a ladder system on $E^{\lambda^+}_{\lambda}$ that enjoys the uniformization property.
Hence, we now turn to discuss the uniformization property at successor of singulars.
We commence with revealing the richer content of Theorem \ref{112}.

\begin{thm}[Shelah, \cite{sh186}] Suppose $\ch_\lambda$ holds for a strong limit singular cardinal, $\lambda$.
If $S\s E^{\lambda^+}_{\cf(\lambda)}$ is a non-reflecting stationary set,
then there exists a notion of forcing $\mathbb P_S$ such that:
\begin{enumerate}
\item $\mathbb P_S$ is $\lambda$-distributive;
\item $\mathbb P_S$ satisfies the $\lambda^{++}$-c.c.;
\item $S$ remains stationary in $V^{\mathbb P_S}$;
\item in $V^{\mathbb P_S}$, there exists a ladder system on $S$ that has the uniformization property.
\end{enumerate}
\end{thm}

By Theorem \ref{124}, it is consistent that diamond fails on a set that reflects stationarily often.
Now, what about the following strengthening:
\begin{q}\index{ZZZ@Open Problems!Question 04} Is it consistent with $\gch$ that for some singular cardinal, $\lambda$,
there exists a stationary set $S\s E^{\lambda^+}_{\cf(\lambda)}$ \emph{that reflects stationarily often},
and a ladder system on $S$ that has the uniformization property?
\end{q}
\begin{remark} By Corollary \ref{121}, $\sap_\lambda$ necessarily fails in such an hypothetical model.\end{remark}

Now, what about the existence of ladder systems that do \emph{not} enjoy the uniformization property?
Clearly, if $\lambda$ is a strong limit singular cardinal, then Theorem \ref{43} does not apply.
For this, consider the following.

\begin{fact}[Shelah, \cite{sh667}] Suppose $\ch_\lambda$ holds for a strong limit singular cardinal, $\lambda$.
Then, for every stationary $S\s\lambda^+$, there exists a ladder system on $S$ that does not enjoy the uniformization property.
\end{fact}
\begin{proof} Fix a stationary $S\s\lambda^+$. If $S\cap E^{\lambda^+}_{\not=\cf(\lambda)}$ is
stationary, then by Theorem \ref{13}, $\diamondsuit_S$ holds, and then by Fact \ref{l25},
moreover, no ladder system on $S$ has the uniformization property.
Next, suppose $S\s E^{\lambda^+}_{\cf(\lambda)}$ is a given stationary set. By the upcoming Theorem \ref{313},
in this case, we may pick a ladder system $\langle L_\alpha\mid \alpha\in S\rangle$ such that
for every function $f:\lambda^+\rightarrow 2$, there exists some $\alpha\in S$
such that if $\{ \alpha_i\mid i<\cf(\lambda)\}$ denotes the increasing enumeration of $L_\alpha$,
then $f(\alpha_{2i})=f(\alpha_{2i+1})$ for all $i<\cf(\lambda)$.

It follows that if for each $\alpha\in S$, we pick $f_\alpha:L_\alpha\rightarrow 2$ satisfying for all $\beta\in L_\alpha$:
$$f_\alpha(\beta)=\begin{cases}0,&\exists i<\cf(\lambda)(\beta=\alpha_{2i})\\1,&\text{otherwise}\end{cases},$$
then the sequence $\langle f_\alpha\mid \alpha\in S\rangle$ cannot be uniformized.
\end{proof}

\begin{remark} Note that the sequence $\langle f_\alpha\mid \alpha\in S\rangle$
that was derived in the preceding proof from the guessing principle of Theorem \ref{313},
is a sequence of non-constant functions that cannot be uniformized.
To compare, the sequence that was derived from weak diamond in the proof of Fact \ref{l25}
is a sequence of constant functions. In other words, weak diamond is stronger in the sense that it entails the existence of
a \emph{monochromatic coloring} that cannot be uniformized.
\end{remark}

\begin{thm}[Shelah, \cite{sh667}]\label{313} Suppose $\ch_\lambda$ holds for a strong limit singular cardinal, $\lambda$,
$S\s E^{\lambda^+}_{\cf(\lambda)}$ is stationary and $\mu<\lambda$ is a given cardinal.

Then there exists a ladder system $\langle L_\alpha\mid \alpha\in S\rangle$
so that if $\{ \alpha_i \mid i<\cf(\lambda)\}$ denotes the increasing enumeration of $L_\alpha$,
then for every function $f:\lambda^+\rightarrow\mu$, the following set is stationary:
$$ \{ \alpha\in S\mid f(\alpha_{2i})=f(\alpha_{2i+1})\text{ for all }i<\cf(\lambda)\}.$$
\end{thm}
\begin{proof} Without loss of generality, $\lambda$ divides the order-type of $\alpha$, for all $\alpha\in S$. Put $\kappa:=\cf(\lambda)$ and $\theta:=2^\kappa$.
By $2^\lambda=\lambda^+$,
let $\{ d_\gamma\mid \gamma<\lambda^+\}$ be some enumeration of $\{ d:\theta\times\tau\rightarrow \mu\mid \tau<\lambda^+\}$.

Fix $\alpha\in S$. Let $\langle c^\alpha_i\mid i<\kappa\rangle$  be the increasing enumeration
of some club subset of $\alpha$, such that $(c^\alpha_i,c^\alpha_{i+1})$ has cardinality $\lambda$ for all $i<\kappa$.
Also, let $\{ b^\alpha_i\mid i<\kappa\}\s[\alpha]^{<\lambda}$ be a continuous chain converging to $\alpha$
with $b^\alpha_i\s c^\alpha_i$
for all $i<\kappa$.
Recall that we have fixed $\alpha\in S$; now, in addition, we also fix $i<\kappa$.

For all $j<\kappa$, define a function $\psi_j=\psi_{\alpha,i,j}: (c^\alpha_i,c^\alpha_{i+1})\rightarrow {}^{\theta\times b^\alpha_j}(\mu+1)$ such that
for all $\varepsilon\in (c^\alpha_i,c^\alpha_{i+1})$ and $(\beta,\gamma)\in\theta\times b^\alpha_j$:
 $$\psi_j(\varepsilon)(\beta,\gamma)=\begin{cases}d_\gamma(\beta,\varepsilon),&(\beta,\varepsilon)\in\dom(d_\gamma)
 \\\mu,&\text{otherwise}\end{cases}.$$
For all $j<\kappa$, since $|{}^{\theta\times b^\alpha_j}(\mu+1)|<\lambda=|(c^\alpha_i,c^\alpha_{i+1})|$ ,
let us pick two ordinals $\alpha^j_{i,0},\alpha^j_{i,1}$ with $c^\alpha_i<\alpha^j_{i,0}<\alpha^j_{i,1}<c^\alpha_{i+1}$
such that $\psi_j(\alpha^j_{i,0})=\psi_j(\alpha^j_{i,1})$.

For every function  $g\in{}^\kappa\kappa$,
consider the ladder system $\langle L_\alpha^g\mid \alpha\in S\rangle$,
where  $L_\alpha^g:=\{ \alpha^{g(i)}_{i,0},\alpha^{g(i)}_{i,1}\mid i<\kappa\}$.

\begin{claim}  There exists some $g\in{}^\kappa\kappa$  such that $\langle L^g_\alpha\mid \alpha\in S\rangle$ works.
\end{claim}
\begin{proof} Suppose not. Let $\{ g_\beta\mid \beta<\theta\}$ be some enumeration of ${}^\kappa\kappa$.
Then, for all $\beta<\theta$, we may pick a function $f_\beta:\lambda^+\rightarrow\mu$
and a club $E_\beta$ such that for all $\alpha\in S\cap E_\beta$, there exists some $i<\kappa$
such that
$$f_\beta(\alpha^{g_\beta(i)}_{i,0})\not=f_\beta(\alpha^{g_\beta(i)}_{i,1}).$$
Now, let $h:\lambda^+\rightarrow\lambda^+$ be the function such that for all $\epsilon<\lambda^+$:
$$h(\epsilon)=\min\{ \gamma<\lambda^+\mid \forall(\beta,\varepsilon)\in\theta\times\epsilon\left(d_\gamma(\beta,\varepsilon)\text{ is defined and equals }f_\beta(\varepsilon)\right)\}.$$
Pick $\alpha\in S\cap\bigcap_{\beta<\theta}E_\beta$ such that $h[\alpha]\s\alpha$.

Then we may define a function $g:\kappa\rightarrow\kappa$ by letting:
$$g(i):=\min\{ j<\kappa\mid h(c^\alpha_{i+1})\in b^\alpha_j\}.$$
Let $\beta<\theta$ be such that $g=g_\beta$ and fix $i<\kappa$ such that
$$f_\beta(\alpha^{g_\beta(i)}_{i,0})\not=f_\beta(\alpha^{g_\beta(i)}_{i,1}).$$

Put $j:=g(i)$. By definition of $\alpha^{j}_{i_,0}$ and $\alpha^{j}_{i_,1}$,
we know that $\psi_{\alpha,i,j}(\alpha^{j}_{i_,0})=\psi_{\alpha,i,j}(\alpha^{j}_{i_,1})$
is a function from $\theta\times b^\alpha_j$ to $\mu+1$.

Put $\gamma:=h(c^\alpha_{i+1})$; then $(\beta,\gamma)\in \theta\times b^\alpha_j$, and hence:
$$\psi_{\alpha,i,j}(\alpha^{j}_{i_,0})(\beta,\gamma)=\psi_{\alpha,i,j}(\alpha^{j}_{i_,1})(\beta,\gamma).$$
It now follows from $\alpha^{j}_{i_,0}<\alpha^{j}_{i_,1}<c^\alpha_{i+1}$ and $\gamma=h(c^\alpha_{i+1})$, that:
$$f_\beta(\alpha^{j}_{i_,0})=d_\gamma(\beta,\alpha^j_{i,0})=\psi_{\alpha,i,j}(\alpha^{j}_{i_,0})(\beta,\gamma)=\psi_{\alpha,i,j}(\alpha^{j}_{i_,1})(\beta,\gamma)=d_\gamma(\beta,\alpha^j_{i,1})=f_\beta(\alpha^{j}_{i_,1})$$
Unrolling the notation, we must conclude that
$$f_\beta(\alpha^{g_\beta(i)}_{i,0})=f_\beta(\alpha^{j}_{i_,0})=f_\beta(\alpha^{j}_{i_,1})=f_\beta(\alpha^{g_\beta(i)}_{i,1}),$$
thus, yielding a contradiction to $\alpha\in E_\beta$.
\end{proof}
Thus, it has been established that there exists a ladder system with the desired properties.
\end{proof}

In light of Theorem \ref{19}, the moral of Theorem \ref{313} is that $\gch$ entails some of the consequences of diamond,
even in the case that  diamond fails.
Two natural questions concerning this theorem are as follows.
\begin{q}\index{ZZZ@Open Problems!Question 05} Is it possible to eliminate the ``strong limit'' hypothesis from Theorem \ref{313},
while maintaining the same conclusion?
\end{q}

\begin{q}\index{ZZZ@Open Problems!Question 06} Is Theorem \ref{313} true also for the case that $\mu=\lambda$?
\end{q}

Note that an affirmative answer to the last question follows from $\diamondsuit_S$.
In fact, even if $2^\lambda>\lambda^+$, but  Ostaszewski's principle, $\clubsuit_S$, holds,
then a ladder system as in Theorem \ref{313} for the case $\mu=\lambda$, exists.
\begin{defn}[Ostaszewski, \cite{ostaszewski}]\label{215} Let $\lambda$ denote an infinite cardinal, and $S$ denote a stationary subset of $\lambda^+$.
Consider the following principle.
\begin{itemize}
\item[$\blacktriangleright$]\index{Guessing Principles!Z2@$\clubsuit_S$} $\clubsuit_S$ asserts that there exists a sequence $\langle  A_\alpha\mid \alpha\in S\rangle$ such that:
\begin{itemize}
\item[$\bullet$] for all $\alpha\in S$, $A_\alpha$ is a cofinal subset of $\alpha$;
\item[$\bullet$] if $Z$ is a \emph{cofinal} subset of $\lambda^+$, then the following set is stationary:
$$\left\{\alpha\in S\mid A_\alpha\s Z\right\}.$$
\end{itemize}
\end{itemize}
\end{defn}

It is worth mentioning that unlike $\clubsuit^-_S$, the principle $\clubsuit_S$ makes sense also in the case that $S\s E^{\lambda^+}_{\lambda}$.
In particular, the missing case of Theorem \ref{12} may be compensated by the observation that
$\diamondsuit_S$ is equivalent to $\clubsuit_S+\ch_\lambda$.
It is also worth mentioning that $\clubsuit_{\lambda^+}+\neg\ch_\lambda$ is indeed consistent; for instance,
in \cite{sh98}, Shelah introduces a model of $\clubsuit_{\omega_1}+\neg\Phi_{\omega_1}$.

Next, consider Theorem \ref{313} for the case that $\mu=\cf(\lambda)$. In this case, the theorem yields
a collection $\mathcal L\s[\lambda^+]^{\cf(\lambda)}$ of size $\lambda^+$, such that for every function
$f:\lambda^+\rightarrow\cf(\lambda)$, there exists some $L\in\mathcal L$ such that $f\restriction L$ is not injective (in some strong sense).
Apparently, this fact led Shelah and D\v{z}amonja  to consider the following dual question.

\begin{qq} Suppose $\lambda$ is a strong limit singular cardinal.

Must there exist a collection $\mathcal P\s[\lambda^+]^{\cf(\lambda)}$
of size $\lambda^+$ such that for every function $f:\lambda^+\rightarrow\cf(\lambda)$
which is non-trivial in the sense that
 $\bigwedge_{\beta<\cf(\lambda)}|f^{-1}\{\beta\}|=\lambda^+$, there exists some $a\in\mathcal P$
such that $f\restriction a$ is injective?
\end{qq}

We shall be concluding this section by describing the resolution of the above question.
To refine the question, consider the following two definitions.
\begin{defn} For a function $f:\lambda^+\rightarrow\cf(\lambda)$, let $\mathfrak \kappa_f$
denote the minimal cardinality of a family $\mathcal P\s[\lambda^+]^{\cf(\lambda)}$
with the property that whenever $Z\s\lambda^+$ satisfies $\bigwedge_{\beta<\cf(\lambda)}|Z\cap f^{-1}\{\beta\}|=\lambda^+$,
then there exist some $a\in\mathcal P$ with $\sup(f[a\cap Z])=\cf(\lambda)$.\footnote{Note that if $\lambda$ is a strong limit, then we may assume that $\mathcal P$ is closed under taking subsets.
Thus, we may moreover demand the existence of  $a\in\mathcal P$ such that $a\s Z$ and $f\restriction a$ is injective.}
\end{defn}

\begin{defn}\index{Guessing Principles!$($@$\lambda^+$-guessing} For a singular cardinal $\lambda$, we say that $\lambda^+$-\emph{guessing}
holds iff $\kappa_f\le\lambda^+$ for all $f\in{}^{\lambda^+}\cf(\lambda)$.
\end{defn}
Answering the above-mentioned question in the negative,
Shelah and D\v{z}amonja  established the consistency of the failure of $\lambda^+$-guessing.
\begin{thm}[D\v{z}amonja-Shelah, \cite{sh685}]\label{ds} It is relatively consistent with the existence of a supercompact cardinal
that there exist a strong limit singular cardinal $\lambda$ and a function $f:\lambda^+\rightarrow\cf(\lambda)$
such that $\kappa_f=2^\lambda>\lambda^+$.
\end{thm}

Recently, we realized that the above-mentioned question is simply equivalent to the question
of whether every strong limit singular cardinal $\lambda$ satisfies $\ch_\lambda$.
\begin{thm}[\cite{gitik-rinot}] Suppose $\lambda$ is a strong limit singular cardinal.
Then: $$\{ \kappa_f\mid f\in{}^{\lambda^+}\cf(\lambda)\}=\{0,2^\lambda\}.$$
\end{thm}

In particular, if $\lambda$ is a strong limit singular cardinal,
then $\lambda^+$-guessing happens to be equivalent to the, seemingly, much stronger principle, $\diamondsuit^+_{E^{\lambda^+}_{\not=\cf(\lambda)}}$.

\newpage

\section{The Souslin Hypothesis and Club Guessing}\label{sectionSH}

Recall that a $\lambda^+$-Aronszajn tree is a tree of height $\lambda^+$,
of width $\lambda$, and without chains of size $\lambda^+$.
A $\lambda^+$-Souslin tree is a $\lambda^+$-Aronszajn tree that has no antichains of size $\lambda^+$.

Jensen introduced the diamond principle and studied its relation to Souslin trees.

\begin{thm}[Jensen, \cite{jensen}]\label{41} If $\lambda^{<\lambda}=\lambda$ is a regular cardinal such that $\diamondsuit_{E^{\lambda^+}_\lambda}$ holds,
then there exists a $\lambda^+$-Souslin tree.

In particular, $\diamondsuit_{\omega_1}$ entails the existence of an $\omega_1$-Souslin tree.
\end{thm}

\begin{thm}[Jensen, see \cite{devlins}]\label{42}  $\gch$ is consistent together with the non-existence of an $\omega_1$-Souslin tree.
\end{thm}

\begin{remark} This is how Jensen proves Theorem \ref{14}. For a more modern proof of Theorem \ref{42}, see  \cite{pid} or \cite{sh81}.
\end{remark}

Let $V$ denote the model from Theorem \ref{14}/\ref{42},
and let $\mathbb P:=\add(\omega,1)$ denote Cohen's notion of forcing for introducing a single Cohen real.
Since $V\models\neg\diamondsuit_{\omega_1}$ and since $\mathbb P$ is c.c.c.,
the discussion after Definition \ref{def13} shows that $V^{\mathbb P}\models\neg\diamondsuit_{\omega_1}$.
By a theorem of Shelah from \cite{sh176}, adding a Cohen real introduces an $\omega_1$-Souslin tree,
and hence $V^{\mathbb P}$ is a model of $\ch$ witnessing the fact that the existence of an $\omega_1$-Souslin tree
does not entail $\diamondsuit_{\omega_1}$.

Now, one may wonder what is the role of the cardinal arithmetic assumption in Theorem \ref{41}?
the answer is that this hypothesis is necessary.
To exemplify the case $\lambda=\aleph_1$,  we mention that $\pfa$ implies $\diamondsuit^+_{E^{\lambda^+}_{\lambda}}$,
but it also implies that $\lambda^{<\lambda}\not=\lambda$ and the \emph{non-existence} of $\lambda^+$-Aronszajn trees.\footnote{For
an introduction to the Proper Forcing Axiom ($\pfa$), see \cite{devlinpfa}.}

So, $\diamondsuit_{E^{\omega_2}_{\omega_1}}$ \emph{per se} does not impose the existence of an $\omega_2$-Souslin tree.
Also, starting with a weakly compact cardinal, Laver and Shelah \cite{sh104}
established that $\ch$ is consistent together with the non-existence of an $\aleph_2$-Souslin tree.
This leads us to the following tenacious question.

\begin{q}[folklore]\label{q20}\index{ZZZ@Open Problems!Question 07} Does $\gch$ imply the existence of an $\omega_2$-Souslin tree?
\end{q}

An even harder question is suggested by Shelah in \cite{sh666}.

\begin{q}[Shelah]\index{ZZZ@Open Problems!Question 08} Is it consistent that the $\gch$ holds while
for some regular uncountable $\lambda$,
there exists neither $\lambda^+$-Souslin trees nor $\lambda^{++}$-Souslin trees?
\end{q}

Gregory's proof of Theorem \ref{15} appears in the paper \cite{gregory} that deals with Question \ref{q20},
and in which this theorem is utilized to supply the following partial answer.

\begin{thm}[Gregory, \cite{gregory}]\label{403} Assume $\gch$ (or just $\ch_\omega+\ch_{\omega_1}$).

If there exists a non-reflecting stationary subset of $E^{\omega_2}_\omega$,
then there exists an $\omega_2$-Souslin tree.
\end{thm}

It follows that the consistency strength of a negative answer to Question \ref{q20} is at least that of the existence of a Mahlo cardinal.
Recently, B. Koenig suggested an approach to show that the strength is at least that of the existence of a weakly compact cardinal.
Let $\square(\omega_2)$\index{Square Principles!$\sap_{z}$@$\square(\omega_2)$} denote the assertion that there exists a sequence $\langle C_\alpha\mid \alpha<\omega_2\rangle$
such that for all limit $\alpha<\omega_2$: (1) $C_\alpha$ is a club subset of $\alpha$,  (2) if $\beta$
is a limit point of $\alpha$, then $C_\alpha\cap\beta=C_\beta$, (3) there exists no ``trivializing'' club $C\s\omega_2$
such that $C\cap\beta=C_\beta$ for all limit points $\beta$ of $C$.

The principle $\square(\omega_2)$ is a consequence of $\square_{\omega_1}$,%
\footnote{\index{Square Principles!$\square_\lambda$}The square property at $\lambda$, denoted $\square_\lambda$,
is the principle $\square_{\lambda,1}$ as in Definition \ref{45}.}
but its consistency strength
is higher --- it is that of the existence of a weakly compact cardinal. Thus, Koenig's question is as follows.

\begin{q}[B. Koenig]\index{ZZZ@Open Problems!Question 09} Does $\gch+\square(\omega_2)$ imply the existence of an $\omega_2$-Souslin tree?
\end{q}

In  light of Theorem \ref{41}, to answer Question \ref{q20} in the affirmative,
one probably needs to find a certain consequence of $\diamondsuit_{E^{\omega_2}_{\omega_1}}$
that, from one hand, follows outright from $\gch$ and which is, on the other hand, strong enough to allow the construction of an $\aleph_2$-Souslin tree.
An example of $\zfc$-provable consequences of diamond is Shelah's family of \emph{club guessing} principles.
The next theorem exemplifies only a few out of many results in this direction.

\begin{thm}[several authors]\label{44}\index{Guessing Principles!ZC@Club guessing} For infinite cardinals $\mu\le\lambda$,
and a stationary set $S\s E^{\lambda^+}_\mu$, there exists a sequence $\overrightarrow{C}=\langle C_\alpha\mid \alpha\in S\rangle$
such that for all $\alpha\in S$, $C_\alpha$ is a club in $\alpha$ of order-type $\mu$, and:
\begin{enumerate}
\item if $\mu<\lambda$, then $\overrightarrow{C}$ may be chosen
such that for every club $D\s\lambda^+$, the following set is stationary:
$$\{ \alpha\in S\mid C_\alpha\s D\}.$$
\item if $\omega<\mu=\cf(\lambda)<\lambda$, then $\overrightarrow{C}$ may be chosen
such that for almost all $\alpha\in S$,
$\langle \cf(\beta)\mid \beta\in\nacc(C_\alpha)\rangle$ is a strictly increasing sequence cofinal in $\lambda$,
 and for every club $D\s\lambda^+$, the following set is stationary:
$$\{ \alpha\in S\mid C_\alpha\s D\}.$$
\item if $V=L$, then $\overrightarrow{C}$ may be chosen such that for every club $D\s\lambda^+$,
the following set contains a club subset of $S$:
$$\{ \alpha\in S\mid \exists\beta<\alpha(C_\alpha\bks\beta\s D)\}.$$
\item if $\omega<\cf(\mu)=\lambda$, then $\overrightarrow{C}$ may be chosen
such that for every club $D\s\lambda^+$, the following set is stationary:
$$\{ \alpha\in S\mid \{\beta\in C_\alpha\mid \min(C_\alpha\bks\beta+1)\in D\}\text{ is stationary in }\alpha \}.$$
\end{enumerate}
\end{thm}

$\blacktriangleright$ Theorem \ref{44}(1) is due to Shelah \cite{shg}, and the principle appearing there reflects the most naive form of club guessing.
Personally, we are curious whether the guessing may concentrate on a prescribed stationary set $T$:
\begin{q}\index{ZZZ@Open Problems!Question 10} Suppose that $S,T$ are given stationary subsets of a successor cardinal $\lambda^+$.
Must there exist a sequence $\langle C_\alpha\mid \alpha\in S\rangle$ with $\sup(C_\alpha)=\alpha$
for all $\alpha\in S$, such that for every club $D\s\lambda^+$, $\{\alpha\in S\mid C_\alpha\s D\cap T\cap\alpha\}$ is stationary?
\end{q}
A positive answer follows from $\clubsuit^-_S$, and a negative answer is consistent for various cardinals $\lambda$ and non-reflecting sets $S\s\lambda^+$,
hence one should focus on sets $S\s E^{\lambda^+}_{\cf(\lambda)}$ that reflect stationarily often, and, e.g., $T=\tr(S)$.

$\blacktriangleright$ Theorem \ref{44}(2) is due to Shelah \cite{shg}, but see also Eisworth and Shelah \cite{sh819}.
Roughly speaking, the principle appearing there requires that, in addition to the naive club guessing,
the non-accumulation points of the local clubs to  be of high cofinality.
An hard open problem is whether their assertion is valid also in the case of countable cofinality.

\begin{q}[Eisworth-Shelah]\index{ZZZ@Open Problems!Question 11}\label{q11} Suppose that $\lambda$ is a singular cardinal of countable cofinality.
Must there exist a ladder system $\langle L_\alpha\mid \alpha\in E^{\lambda^+}_{\cf(\lambda)}\rangle$
such that for almost all $\alpha$, $\langle \cf(\beta)\mid \beta\in L_\alpha\rangle$ is a strictly increasing
$\omega$-sequence cofinal in $\lambda$, and for every club $D\s\lambda^+$,
the set $\{\alpha\in E^{\lambda^+}_{\cf(\lambda)}\mid L_\alpha\s D\}$
is stationary?
\end{q}

While the above question remains open, Eisworth recently established the validity of a principle named \emph{off-center club guessing} \cite{eisworthapal},
and demonstrated that the new principle can serve as a useful substitute to the principle of Question \ref{q11}.

$\blacktriangleright$ Theorem \ref{44}(3) is due to Ishiu \cite{MR2194236},
and the principle appearing there is named \emph{strong club guessing}\index{Guessing Principles!ZSC@Strong club guessing}.
The ``strong'' stands for the requirement that the guessing is done on almost all points rather that on just stationary many.
Historically, Foreman and Komj\'ath first proved in \cite{foremankomjath}
that strong club guessing may be introduced by forcing (See Theorem \ref{418} below),
and later on, Ishiu proved that this follows from $V=L$.
In his paper, Ishiu asks whether $V=L$ may be reduced to a diamond-type hypothesis. Here is a variant of his question.

\begin{q}\index{ZZZ@Open Problems!Question 12} Suppose that $\diamondsuit^+_{\lambda^+}$ holds
for a given infinite cardinal $\lambda$.
Must there exist a regular cardinal $\mu<\lambda$, a stationary set $S\s E^{\lambda^+}_\mu$,
and a ladder system $\langle L_\alpha\mid \alpha\in S\rangle$ such
that for every club $D\s\lambda^+$, for club many $\alpha\in S$, there exists $\beta<\alpha$ with $L_\alpha\bks\beta\s D$?
\end{q}

We mention that $\diamondsuit^+_{\omega_1}$ is consistent together with
the failure of strong club guessing over $\omega_1$ (see \cite{plusiteration}), while, for an uncountable regular cardinal $\lambda$,
and a stationary $S\s E^{\lambda^+}_\lambda$, $\diamondsuit^*_S$ suffices to yield strong club guessing over $S$.

$\blacktriangleright$ Theorem \ref{44}(4) is due to Shelah \cite{sh572},
and a nice presentation of the proof may be found in \cite{soukups}.
The prototype of this principle is
the existence of a sequence of local clubs, $\langle C_\alpha\mid \alpha\in E^{\lambda^+}_{\lambda}\rangle$,
such that for every club $D\s\lambda^+$,
there exists some $\alpha\in E^{\lambda^+}_\lambda$ with $\sup(\nacc(C_\alpha)\cap D)=\alpha$.
Now, if $\{ \alpha_i \mid i<\lambda\}$ denotes the increasing enumeration of $C_\alpha$,
then Theorem \ref{44}(3) states that for every club $D\s\lambda^+$, there exists stationarily many $\alpha\in S$,
for which not only that $\sup(\nacc(C_\alpha)\cap D)=\alpha$, but moreover,
$\{ i<\lambda\mid \alpha_{i+1}\in D\}$ is stationary in $\lambda$.
According to Shelah \cite{sh666}, to answer Question \ref{q20} in the affirmative,
it suffices to find a proof of the following natural improvement.
\begin{q}[Shelah]\index{ZZZ@Open Problems!Question 13} For a regular uncountable cardinal, $\lambda$, must there exist
a sequence $\langle C_\alpha\mid \alpha\in E^{\lambda^+}_{\lambda}\rangle$
with each $C_\alpha$ a club in $\alpha$
whose increasing enumeration is $\{ \alpha_i\mid i<\lambda\}$,
such that for every club $D\s\lambda^+$,
there exists stationarily many $\alpha$,
for which $\{ i<\lambda\mid \alpha_{i+1}\in D\text{ and }\alpha_{i+2}\in D\}$ is stationary in $\lambda$?
\end{q}

To exemplify the tight relation between higher Souslin trees and the preceding type of club guessing,
we mention the next principle.
\begin{defn}[\cite{rinot09}]\index{Guessing Principles!$(_)$@$\langle T\rangle_S$} Suppose  $\lambda$ is a regular uncountable cardinal, $T$ is a stationary subset of $\lambda$,
and $S$ is a stationary subset of $E^{\lambda^+}_\lambda$.

$\langle T\rangle_S$ asserts the existence of 
sequences $\langle C_\alpha\mid \alpha\in S\rangle$ and $\langle A^\alpha_i\mid \alpha\in S, i<\lambda\rangle$
such that:
\begin{enumerate}
\item for all $\alpha\in S$, $C_\alpha$ is  a club subset of $\alpha$ of order-type $\lambda$;
\item if for all $\alpha\in S$, $\{ \alpha_i\mid i<\lambda\}$ denotes the increasing enumeration of $C_\alpha$,
then for every club $D\s\lambda^+$ and every subset $A\s\lambda^+$, there exist stationarily many $\alpha\in S$ for which:
$$\{ i\in T\mid \alpha_{i+1}\in D\ \&\ A\cap\alpha_{i+1}=A^{\alpha}_{i+1}\}\text{ is stationary in }\lambda.$$
\end{enumerate}
\end{defn}

It is obvious that $\diamondsuit_S\gorer\langle T\rangle_S$.
It is also not hard to see that
$\langle T\rangle_S\gorer\diamondsuit_S$ whenever $\ns_\lambda\restriction T$ is saturated.\footnote{See Definition \ref{40} below.}
A strengthening of Theorem \ref{41}  is the following.
\begin{thm}[implicit in \cite{sh449}] If $\lambda^{<\lambda}=\lambda$ is a regular uncountable cardinal and $\langle\lambda\rangle_{E^{\lambda^+}_\lambda}$ holds,
then there exists a $\lambda^+$-Souslin tree.
\end{thm}

We now turn to discuss Souslin trees at the of successor of singulars.
By Magidor and Shelah \cite{sh324}, if $\lambda$ is a singular cardinal
which is a limit of strongly compact cardinals, then there are no $\lambda^+$-Aronszajn trees.
In particular, it is consistent with $\gch$ that for some singular cardinal $\lambda$,
there are no $\lambda^+$-Souslin trees.
On the other hand, Jensen proved the following.
\begin{thm}[Jensen]\label{35} For a singular cardinal $\lambda$, $\ch_\lambda+\square_\lambda$ entails the existence of
a $\lambda^+$-Souslin tree.
\end{thm}

Since $\square_\lambda\gorer\square^*_\lambda$ and the latter still witnesses the existence of a $\lambda^+$-Aronszajn tree,
the question which appears to be the agreed analogue of Question \ref{q20} is the following.

\begin{q}[folklore]\label{q19}\index{ZZZ@Open Problems!Question 14} For a singular cardinal $\lambda$, does $\gch+\square^*_\lambda$ imply the existence of a $\lambda^+$-Souslin tree?
\end{q}

A minor modification of Jensen's proof of Theorem \ref{35} entails
a positive answer to Question \ref{q19} provided that there exists a non-reflecting stationary subset of $E^{\lambda^+}_{\not=\cf(\lambda)}$.
However, by Magidor and Ben-David \cite{magidorbendavid}, it is relatively consistent with the
existence of a supercompact cardinal that the $\gch$ holds, $\square^*_{\aleph_\omega}$ holds,
and every stationary subset of $E^{\aleph_{\omega+1}}_{\not=\omega}$ reflects.

A few years ago, Schimmerling \cite{schsous} suggested that the
community should perhaps try to attack a softer version of Question \ref{q19},
which is related to the following hierarchy of square principles.

\begin{defn}[Schimmerling, \cite{schimmerling}]\label{45}\index{Square Principles!$\square_{\lambda,<\mu}$} For cardinals, $\mu,\lambda$, $\square_{\lambda,<\mu}$
asserts the existence of a sequence $\langle \mathcal C_\alpha\mid \alpha<\lambda^+\rangle$ such that
for all limit $\alpha<\lambda^+$:
\begin{itemize}
\item $0<|\mathcal C_\alpha|<\mu$;
\item $C$ is a club subset of $\alpha$ for all $C\in\mathcal C_\alpha$;
\item if $\cf(\alpha)<\lambda$, then $\otp(C)<\lambda$ for all $C\in\mathcal C_\alpha$;
\item if $C\in\mathcal C_\alpha$ and $\beta\in\acc(C)$, then $C\cap\beta\in\mathcal C_\beta$.
\end{itemize}
We also write $\square_{\lambda,\mu}$ for $\square_{\lambda,<\mu^+}$.
\end{defn}

\begin{q}[Schimmerling]\label{q18}\index{ZZZ@Open Problems!Question 15} Does $\gch+\square_{\aleph_\omega,\omega}$ imply the existence of an $\aleph_{\omega+1}$-Souslin tree?
\end{q}

In \cite{sh221}, Abraham, Shelah and Solovay showed that if $\ch_\lambda+\square_\lambda$
holds for a given strong limit singular cardinal, $\lambda$,
then a principle which is called \emph{square with built-in diamond} may be inferred.
Then, they continued to show how to construct a $\lambda^+$-Souslin tree with a certain special property, based on this principle.

There are several variations of square-with-built-in-diamond principles (the first instance appearing in \cite{gray}),
and several constructions of peculiar trees that utilizes principles of this flavor (see \cite{sh236}, \cite{cummingstree}, \cite{sh449}, \cite{novak}).
Recalling the work of Abraham-Shelah-Solovay in \cite{sh221}, it seems reasonable to seek for a principle
that ramifies the hypothesis of Question \ref{q18}. Here is our humble suggestion.

\begin{defn}[\cite{rinot10}]\label{46}\index{Square Principles!Z@$\ssd_{\lambda,<\mu}$} For cardinals, $\mu,\lambda$, $\sd_{\lambda,<\mu}$
asserts the existence of two sequences, $\langle \mathcal C_\alpha\mid \alpha<\lambda^+\rangle$
and $\langle \varphi_\theta\mid \theta\in \Gamma\rangle$,
such that all of the following holds:
\begin{itemize}
\item $\emptyset\not=\Gamma\s\{\theta<\lambda^+\mid \cf(\theta)=\theta\}$;
\item $\langle \mathcal C_\alpha\mid \alpha<\lambda^+\rangle$ is a $\square_{\lambda,<\mu}$-sequence;
\item $\varphi_\theta:\mathcal P(\lambda^+)\rightarrow\mathcal P(\lambda^+)$ is a function, for all $\theta\in\Gamma$;
\item for every subset $A\s\lambda^+$, every club $D\s\lambda^+$,
and every cardinal $\theta\in\Gamma$,
there exists some $\alpha\in E^{\lambda^+}_\theta$ such that for all $C\in\mathcal C_\alpha$:
$$\sup\{ \beta\in\nacc(\acc(C))\cap D\mid \varphi_\theta(C\cap\beta)=A\cap\beta\}=\alpha.$$
\end{itemize}
We write $\sd_{\lambda,\mu}$ for $\sd_{\lambda,<\mu^+}$.
\end{defn}

Notice that the above principle combines square, diamond and club guessing. The value of this definition is witnessed by the following.

\begin{thm}[\cite{rinot10}]\label{47} Suppose that $\lambda$ is an uncountable cardinal.

If $\sd_{\lambda,\lambda}$ holds, then there exists a $\lambda^+$-Souslin tree.
\end{thm}
\begin{remark}An interesting feature of the (easy) proof of the preceding theorem
is that the construction does not depend on whether $\lambda$ is a regular cardinal or a singular one.
\end{remark}

It follows that if $\gch+\square_{\aleph_\omega,\omega}$
entails $\sd_{\aleph_\omega,\aleph_{\omega}}$,
then this would supply an affirmative answer to Question \ref{q18}.
However, so far, a ramification is available only for the case $\mu\le\cf(\lambda)$.

\begin{thm}[\cite{rinot10}]\label{48} For cardinals $\lambda\ge\aleph_2$, and $\mu\le\cf(\lambda)$, the following are equivalent:
\begin{itemize}
\item[(a)] $\square_{\lambda,<\mu}+\ch_\lambda$;
\item[(b)] $\sd_{\lambda,<\mu}$.
\end{itemize}
\end{thm}
\begin{remark} In the proof of (a)$\gorer$(b), we obtain a $\sd_{\lambda,<\mu}$-sequence
as in Definition \ref{46}
for which, moreover, $\Gamma$ is a non-empty \emph{final segment} of $\{\theta<\lambda\mid \cf(\theta)=\theta\}$.
\end{remark}

Clearly, in the presence of a non-reflecting stationary set, one can push Theorem \ref{48} much further (Cf. \cite{rinot10}).
Thus, to see the difficulty of dealing with the case $\mu=\cf(\lambda)^+$, consider the following variation of club guessing.

\begin{q}\index{ZZZ@Open Problems!Question 16} Suppose that $\lambda$ is a singular cardinal, $\square_{\lambda,\cf(\lambda)}$ holds,
and every stationary subset of $\lambda^+$ reflects.

Must there exist a regular cardinal $\theta$ with $\cf(\lambda)<\theta<\lambda$
 and a $\square_{\lambda,\cf(\lambda)}$-sequence, $\langle \mathcal C_\alpha\mid \alpha<\lambda^+\rangle$,
 such that for every club $D\s\lambda^+$, there exists some $\alpha\in E^{\lambda^+}_\theta$
satisfying $\sup(\nacc(C)\cap D)=\alpha$ for all $C\in\mathcal C_\alpha$?
\end{q}

To conclude this section, let us mention two questions that suggests an alternative generalizations of Theorems \ref{41} and \ref{35}.

\begin{q}[Juh\'asz]\index{ZZZ@Open Problems!Question 17} Does $\clubsuit_{\omega_1}$ entail the existence of an $\aleph_1$-Souslin tree?
\end{q}

\begin{q}[Magidor]\index{ZZZ@Open Problems!Question 18} For a singular cardinal $\lambda$, does $\square_\lambda$
entail the existence of a $\lambda^+$-Souslin tree?
\end{q}

Juh\'asz's question is well-known and a description of its surrounding results deserves a survey paper of its own.
Here, we just mention that most of these results may be formulated in terms of the parameterized diamond principles of \cite{MR2048518}.
For instance, see \cite{heike}.

To answer Magidor's question, one needs to find a yet another $\gch$-free version of diamond
which suggests some non-trivial guessing features.
In \cite{sh775},
Shelah introduced  a principle of this flavor, named \emph{Middle Diamond}\index{Guessing Principles!ZZ@Middle Diamond},
and a corollary to the results of \cite[$\S4$]{sh829} reads as follows (compare with Definitions \ref{0} and \ref{215}.)

\begin{thm}[Shelah, \cite{sh829}]For every cardinal $\lambda\ge\beth_{\omega_1}$,
there exist a finite set $\mathfrak d\s\beth_{\omega_1}$,
and a sequence $\langle (C_\alpha,A_\alpha) \mid \alpha<\lambda^+\rangle$ such that:
\begin{itemize}
\item for all limit $\alpha$,  $C_\alpha$ is a club in $\alpha$, and $A_\alpha\s C_\alpha$;
\item[$\bullet$] if $Z$ is a subset of $\lambda^+$,
then for every regular cardinal $\kappa\in\beth_{\omega_1}\bks\mathfrak d$,
the following set is stationary:
$$\{\alpha\in E^{\lambda^+}_\kappa \mid Z\cap C_\alpha=A_\alpha\}.$$
\end{itemize}
\end{thm}
For more information on the middle diamond, consult \cite{rinotr01}.

\section{Saturation of the Nonstationary Ideal}\label{sectionaturation}

\begin{defn}[folklore]\label{40}\index{Anti-$\diamondsuit$ Principles!$\ns_{\lambda^+}\restriction S$ is saturated}Suppose that $S$ is a stationary subset of a cardinal, $\lambda^+$.
We say that $\ns_{\lambda^+}\restriction S$ is \emph{saturated} iff for any family $\mathcal F$ of $\lambda^{++}$
many stationary subsets of $S$, there exists two distinct sets $S_1,S_2\in\mathcal F$
such that $S_1\cap S_2$ is stationary.

Of course, we say that $\ns_{\lambda^+}$ is saturated iff $\ns_{\lambda^+}\restriction\lambda^+$ is saturated.
\end{defn}

Now, suppose that $\diamondsuit_S$ holds, as witnessed by $\langle A_\alpha\mid \alpha\in S\rangle$.
For every subset $Z\s\lambda^+$, consider the set $\mathcal G_Z:=\{\alpha\in S\mid Z\cap\alpha=A_\alpha\}$.
Then $\mathcal G_Z$ is stationary and $|\mathcal G_{Z_1}\cap\mathcal G_{Z_2}|<\lambda^+$ for all distinct $Z_1,Z_2\in\mathcal P(\lambda^+)$.
Thus, $\diamondsuit_S$ entails that $\ns_{\lambda^+}\restriction S$ is non-saturated.
For stationary subsets of $E^{\lambda^+}_{<\lambda}$,
an indirect proof of this last observation follows from Theorem \ref{22} below.
For this, we first remind our reader that a set $\mathcal X\s\mathcal P(\lambda^+)$
is said to be \emph{stationary} (in the generalized sense) iff for any function $f:[\lambda^+]^{<\omega}\rightarrow\lambda^+$
there exists some $X\in\mathcal X$ with $f``[X]^{<\omega}\s X$.

\begin{defn}[Gitik-Rinot, \cite{gitik-rinot}]  For an infinite cardinal $\lambda$ and a stationary set $S\s\lambda^+$,
consider the following two principles.

\begin{itemize}

\item[$\blacktriangleright$]\index{Guessing Principles!$(1)_S$} $(1)_S$ asserts that there exists a stationary $\mathcal X\s[\lambda^+]^{<\lambda}$ such that:
\begin{itemize}
\item[$\bullet$] the $\sup$-function on $\mathcal X$ is 1-to-1;
\item[$\bullet$] $\{ \sup(X)\mid X\in\mathcal X\}\s S$.
\end{itemize}

\item[$\blacktriangleright$]\index{Guessing Principles!$(\lambda)_S$}  $(\lambda)_S$ asserts that there exists a stationary $\mathcal X\s[\lambda^+]^{<\lambda}$ such that:
\begin{itemize}
\item[$\bullet$] the $\sup$-function on $\mathcal X$ is $(\le\lambda)$-to-1;
\item[$\bullet$] $\{ \sup(X)\mid X\in\mathcal X\}\s S$.
\end{itemize}
\end{itemize}
\end{defn}

\begin{thm}\label{22} For an uncountable cardinal $\lambda$, and a stationary set $S\s E^{\lambda^+}_{<\lambda}$,
the implication $(1)\gorer(2)\gorer(3)\gorer(4)\gorer(5)$ holds:
\begin{enumerate}
\item $\diamondsuit_S$;
\item $(1)_S$;
\item $(\lambda)_S$;
\item $\clubsuit^-_S$;
\item $\ns_{\lambda^+}\restriction S$ is non-saturated.
\end{enumerate}
\end{thm}
\begin{proof} For a proof of the implication $(1)\gorer(2)\gorer(3)\gorer(4)$, see \cite{gitik-rinot}.
The proof of the last implication appears in \cite{rinot07}, building on the arguments of \cite{sh545}.\end{proof}
Note that by Theorem  \ref{12}, the first four items of the preceding theorem coincide assuming $\ch_\lambda$.
In particular, the next question happens to be the contrapositive version of Question \ref{q1}.
\begin{q}\label{q3}\index{ZZZ@Open Problems!Question 19} Suppose that $\lambda$ is a singular cardinal.
Does $\ch_\lambda$ entail the existence of a stationary $\mathcal X\s [\lambda^+]^{<\lambda}$
on which $X\mapsto\sup(X)$ is an injective map from  $\mathcal X$ to $E^{\lambda^+}_{\cf(\lambda)}$?
\end{q}

Back to non-saturation, since $\zfc\vdash \clubsuit^-_S$ for every stationary subset $S\s E^{\lambda^+}_{\not=\cf(\lambda)}$,
one obtains the following  analogue   of Theorem \ref{13}.

\begin{cor}[Shelah, \cite{shg}] If $\lambda$ is an uncountable cardinal,
and $S$ is a stationary subset of $E^{\lambda^+}_{\not=\cf(\lambda)}$, then $\ns_{\lambda^+}\restriction S$ is non-saturated.
\end{cor}

Thus, as in diamond, we are led to focus our attention on the saturation of $\ns_{\lambda^+}\restriction S$
for stationary sets $S$ which concentrates on the set of critical cofinality.

Kunen \cite{kunensaturated} was the first to establish the consistency of an abstract saturated ideal on $\omega_1$.
As for the saturation of the ideal $\ns_{\omega_1}$, this has been obtained first by Steel and Van Wesep
by forcing over a model of determinacy.

\begin{thm}[Steel-Van Wesep, \cite{adsat}]\label{25} Suppose that $V$ is a model of ``$\zf+\ad_{\mathbb R}+\Theta$ is regular''.
Then, there is a forcing extension satisfying $\zfc+\ns_{\omega_1}$ is saturated.
\end{thm}

Woodin \cite{woodincabal} obtained the same conclusion while weakening the hypothesis to the assumption ``$V=L(\mathbb R)+\ad$''.
Several years later, in \cite{mm1}, Foreman, Magidor and Shelah introduced \emph{Martin's Maximum}, $\mm$,
established its consistency from a supercompact cardinal, and proved that $\mm$
entails that $\ns_{\omega_1}$ is saturated, and remains as such in any $c.c.c.$ extension of the universe.

Then, in \cite{sh253}, Shelah established
the consistency of the saturation of $\ns_{\omega_1}$ from just a Woodin cardinal.
Finally, recent work of Jensen and Steel on the existence of the core model below a Woodin cardinal
yields the following definite resolution.

\begin{thm}[Shelah, Jensen-Steel] The following are equiconsistent:
\begin{enumerate}
\item $\zfc+$``there exist a woodin cardinal'';
\item $\zfc+$``$\ns_{\omega_1}$ is saturated''.
\end{enumerate}
\end{thm}

However, none of these results serves as
a complete analogue of Theorem \ref{14} in the sense
that the following is still open.

\begin{q}[folklore]\label{q6}\index{ZZZ@Open Problems!Question 20} Is $\ch$ consistent with $\ns_{\omega_1}$ being saturated?
\end{q}
\begin{remark} By \cite{rinot10}, ``$\ch+\ns_{\omega_1}$ is saturated'' entails $\sd_{\omega_1,\omega_1}$.
\end{remark}

Recalling that $\ch\gorer\Phi_{\omega_1}$,
it is worth pointing out
that while the saturation of $\ns_{\omega_1}$ is indeed an anti-$\diamondsuit_{\omega_1}$ principle, it is not an anti-$\Phi_{\omega_1}$ principle.
To exemplify this, start with a model of $\mm$ and add $\aleph_{\omega_1}$ many Cohen reals over this model;
then as a consequence of Theorem \ref{311} and the fact that Cohen forcing is $c.c.c.$,
one obtains a model in which $\ns_{\omega_1}$ is still saturated,
while $\Phi_{\omega_1}$ holds.

Let us consider a strengthening of saturation which does serve as an anti-$\Phi_{\lambda^+}$ principle.

\begin{defn}[folklore]\index{Anti-$\diamondsuit$ Principles!$\ns_{\lambda^+}\restriction S$ is dense}Suppose that $S$ is a stationary subset of a cardinal, $\lambda^+$.
We say that $\ns_{\lambda^+}\restriction S$ is \emph{dense} iff there exists a family $\mathcal F$ of $\lambda^{+}$ many stationary
subsets of $S$, such that for any stationary subset $S_1\s S$, there exists some $S_2\in\mathcal F$
such that $S_2\bks S_1$ is non-stationary.

Of course, we say that $\ns_{\lambda^+}$ is dense iff $\ns_{\lambda^+}\restriction\lambda^+$ is dense.
\end{defn}

It is not hard to see that if $\ns_{\lambda^+}\restriction S$ is dense, then it is also saturated.
The above discussion and the next theorem entails that these principles do not coincide.
\begin{thm}[Shelah, \cite{shelaharound}] If $\Phi_{\omega_1}$ holds, then $\ns_{\omega_1}$ is not dense.
\end{thm}

Improving Theorem \ref{25}, Woodin proved:
\begin{thm}[Woodin, \cite{woodinbook}] Suppose that $V$ is a model of $``V=L(\mathbb R)+\ad$''.
Then there is a forcing extension of $\zfc$ in which $\ns_{\omega_1}$ is dense.
\end{thm}

The best approximation for a positive answer to Question \ref{q6} is, as well, due to Woodin,
who proved that $\ch$ is consistent together with $\ns_{\omega_1}\restriction S$ being dense
for some stationary $S\s\omega_1$. Woodin
also obtained an approximation for a negative answer to the very same question.
By \cite{woodinbook}, if $\ns_{\omega_1}$ is saturated and there exists a measurable cardinal, then $\ch$ must fail.

As for an analogue of Theorem \ref{109} --- the following is completely open:

\begin{q}[folklore]\index{ZZZ@Open Problems!Question 21} Is it consistent that $\ns_{\omega_2}\restriction E^{\omega_2}_{\omega_1}$ is saturated?
\end{q}

A major, related, result is the following unpublished theorem of Woodin (for a proof, see \cite[$\S8$]{foremanchapter}.)
\begin{thm}[Woodin]\label{1234} Suppose that $\lambda$ is an uncountable regular cardinal and $\kappa$ is a huge cardinal above it.
Then there exists a $<\lambda$-closed notion of forcing $\mathbb P$, such that in $V^{\mathbb P}$ the following holds:
\begin{enumerate}
\item $\kappa=\lambda^+$;
\item there exists a stationary $S\s E^{\lambda^+}_{\lambda}$ such that $\ns_{\lambda^+}\restriction S$ is saturated.
\end{enumerate}

Moreover, if $\gch$ holds in the ground model, then $\gch$ holds in the extension.
\end{thm}

Foreman, elaborating on Woodin's proof, established the consistency of the saturation of $\ns_{\lambda^+}\restriction S$
for some stationary set $S\s E^{\lambda^+}_\lambda$ and a supercompact cardinal, $\lambda$,
and showed that it is then possible to collapse $\lambda$ to $\aleph_\omega$,
while preserving saturation. Thus, yielding:

\begin{thm}[Foreman] It is relatively consistent with the existence of a supercompact cardinal and an almost huge cardinal above it,
that the $\gch$ holds, and $\ns_{\aleph_{\omega+1}}\restriction S$ is saturated for some stationary $S\s E^{\aleph_{\omega+1}}_\omega$.
\end{thm}

Since the stationary set $S$ was originally a subset of $E^{\lambda^+}_{\lambda}$,
it is a non-reflecting stationary set. This raises the following question.

\begin{q}[folklore]\label{q7}\index{ZZZ@Open Problems!Question 22} Suppose that $\lambda$ is a singular cardinal, and $S\s E^{\lambda^+}_{\cf(\lambda)}$ reflects stationarily often,
must $\ns_{\lambda^+}\restriction S$  be non-saturated?
\end{q}

Recently, the author \cite{rinot07} found several partial answers to Question \ref{q7}.
To start with, as a consequence of Theorem \ref{18} and Theorem \ref{22}, we have:
\begin{thm}[\cite{rinot07}] Suppose $S\s \lambda^+$ is a stationary set,
for a singular cardinal $\lambda$. If $I[S;\lambda]$ contains a stationary set,
then $\ns_{\lambda^+}\restriction S$ is non-saturated.
\end{thm}

In particular, $\sap_\lambda$ (and hence $\square^*_\lambda$) impose a positive answer to Question \ref{q7}.

Recalling Theorem \ref{130}, we also obtain the following.
\begin{thm}[\cite{rinot07}] If $\lambda$ is a singular cardinal of uncountable cofinality and $S\s\lambda^+$
is a stationary set such that $\ns_{\lambda^+}\restriction S$ is saturated,
then for every regular cardinal $\theta$ with $\cf(\lambda)<\theta<\lambda$, at least one of the two holds:
\begin{enumerate}
\item $R_2(\theta,\cf(\lambda))$ fails;
\item $\tr(S)\cap E^{\lambda^+}_\theta$ is nonstationary.
\end{enumerate}
\end{thm}

Next, to describe an additional aspect of Question \ref{q7},
we remind our reader that a set $T\s\lambda^+$ is said to \emph{carry a weak square sequence} iff there
exists sequence $\langle C_\alpha\mid \alpha\in T\rangle$ such that:
\begin{enumerate}
\item $C_\alpha$ is a club subset of $\alpha$ of order-type $\le\lambda$, for all limit $\alpha\in T$;
\item  $|\{ C_\alpha\cap\gamma \mid \alpha\in T\}|\le\lambda$ for all $\gamma<\lambda^+$.
\end{enumerate}

\begin{fact}[\cite{rinot07}]\label{131} Suppose $\lambda$ is a singular cardinal, and $S\s \lambda^+$ is a given stationary set.
If some stationary subset of $\tr(S)$ carries a weak square sequence, then $I[S;\lambda]$ contains a stationary set,
and in particular, $\ns_{\lambda^+}\restriction S$ is non-saturated.
\end{fact}

The consistency of the existence of a stationary set that does not carry a weak square sequence is well-known,
and goes back to Shelah's paper \cite{sh108}. However, the following question is still open.

\begin{q}\index{ZZZ@Open Problems!Question 23} Suppose that $\lambda$ is a singular cardinal. Must there exist a stationary subset of $E^{\lambda^+}_{>\cf(\lambda)}$
that carries a partial weak square sequence?
\end{q}
\begin{remark} The last question is closely related to a conjecture of Foreman and Todorcevic from \cite[$\S6$]{foremantodorcevic}.
Note that by Fact \ref{131} and Theorems \ref{12}, \ref{18}, a positive answer imposes a negative answer on Question \ref{q1}.
\end{remark}

Back to Question \ref{q7}, still, there are a few $\zfc$ results; the first being:

\begin{thm}[Gitik-Shelah, \cite{sh577}] If $\lambda$ is a singular cardinal, then $\ns_{\lambda^+}\restriction E^{\lambda^+}_{\cf(\lambda)}$ is non-saturated.
\end{thm}

Gitik and Shelah's  proof utilizes the $\zfc$ fact that a certain weakening of the club guessing principle from Theorem \ref{44}(2) holds
for all singular cardinal, $\lambda$.
Then, they show that if $\ns_{\lambda^+}\restriction E^{\lambda^+}_{\cf(\lambda)}$
were saturated, then their club guessing principle may be strengthened to a principle that combines their variation of \ref{44}(2),
together with \ref{44}(3). However, as they show, this strong combination is already inconsistent.

In \cite{fatsets}, Krueger pushed further the above argument, yielding the following generalization.
\begin{thm}[Krueger, \cite{fatsets}] If $\lambda$ is a singular cardinal and $S\s\lambda^+$
is a stationary set such that $\ns_{\lambda^+}\restriction S$ is saturated, then $S$ is co-fat.\footnote{Here,
a set $T\s\lambda^+$ is \emph{fat} iff for every cardinal $\kappa<\lambda$ and every club $D\s\lambda^+$,
$T\cap D$ contains some closed subset of order-type $\kappa$.}
\end{thm}

To conclude this section, we mention two complementary results to the Gitik-Shelah argument.

\begin{thm}[Foreman-Komj\'ath, \cite{foremankomjath}]\label{418}  Suppose that $\lambda$ is an uncountable regular cardinal
and $\kappa$ is an almost huge cardinal above it.
Then there exists a notion of forcing $\mathbb P$, such that in $V^{\mathbb P}$ the following holds:
\begin{enumerate}
\item $\kappa=\lambda^+$;
\item there exists a stationary $S\s E^{\lambda^+}_{\lambda}$ such that $\ns_{\lambda^+}\restriction S$ is saturated;
\item $E^{\lambda^+}_\mu$ carries a strong club guessing sequence for any regular $\mu\le\lambda$.
\end{enumerate}
\end{thm}
\begin{remark} By \emph{strong club guessing}, we refer to the principle appearing in Theorem \ref{44}(3).\end{remark}

\begin{thm}[Woodin, \cite{woodinbook}] Assuming $\ad^{L(\mathbb R)}$,  there exists a forcing
extension of $L(\mathbb R)$ in which:
\begin{enumerate}
\item $\ns_{\omega_1}$ is saturated;
\item there exists a strong club guessing sequence on $E^{\omega_1}_\omega$.
\end{enumerate}
\end{thm}
For interesting variations of Woodin's theorem, we refer the reader to \cite{club_pmax}.

\printindex\label{indexpage}
\bibliographystyle{plain}

\end{document}